\documentclass[10pt,leqno,twoside]{article}
\usepackage[a5paper,layoutwidth=148.5mm,layoutheight=7.77817in,outer=10mm,inner=38pt,tmargin=10mm,bmargin=10mm,footskip=15pt]{geometry}
\usepackage[british]{babel}
\usepackage{lmodern}
\usepackage{fix-cm}
\usepackage[T1]{fontenc}
\usepackage{url}
\usepackage[utf8]{inputenc}
\usepackage{xspace}

\let\emph\textsl
\let\textit\textsl

\usepackage{natbib}

\makeatletter
\DeclareFontEncoding{LS1}{}{}
\DeclareFontEncoding{LS2}{}{\noaccents@}
\DeclareSymbolFont{stix-letters}       {LS1}{stix}     {m}{it}
\DeclareSymbolFont{stix-arrows1}       {LS1}{stixsf}   {m} {n}
\DeclareSymbolFont{stix-operators}     {LS1}{stix}     {m} {n}
\DeclareSymbolFont{stix-largesymbols}  {LS2}{stixex}   {m} {n}
\DeclareSymbolFont{stix-bold-operators}{LS1}{stix}     {b} {n}
\DeclareFontSubstitution{LS1}{stix}{m}{n}
\DeclareFontSubstitution{LS2}{stix}{m}{n}
\SetSymbolFont{stix-letters}     {bold}{LS1}{stix}     {b}{it}
\SetSymbolFont{stix-arrows1}     {bold}{LS1}{stixsf}   {b} {n}
\SetSymbolFont{stix-operators}   {bold}{LS1}{stix}     {b} {n}
\SetSymbolFont{stix-largesymbols}{bold}{LS2}{stixex}   {b} {n}

\def\stix@undefine#1{%
    \if\relax\noexpand#1\let#1=\@undefined\fi}
\def\stix@MathSymbol#1#2#3#4{%
    \stix@undefine#1%
    \DeclareMathSymbol{#1}{#2}{#3}{#4}}
\stix@MathSymbol{\stixwedge}{\mathbin}  {stix-operators}{"E1} \let\land=\wedge
\stix@MathSymbol{\stixvee}  {\mathbin}  {stix-operators}{"E2} \let\lor=\vee
\def\bigwedgegras{\DOTSI\bigwedgeop\slimits@}
\def\bigveegras{\DOTSI\bigveeop\slimits@}
\stix@MathSymbol{\stixbigwedgeop}{\mathop}{stix-largesymbols}{"B4}
\stix@MathSymbol{\stixbigveeop}  {\mathop}{stix-largesymbols}{"B5}
\stix@MathSymbol{\stixrightarrow}               {\mathrel}{stix-arrows1}{"99}
\stix@MathSymbol{\stixrightleftarrows}          {\mathrel}{stix-arrows1}{"CB}
\makeatother

\renewcommand\land{\mathbin{\boldsymbol{\stixwedge}}}
\renewcommand\lor{\mathbin{\boldsymbol{\stixvee}}}

\usepackage{amsmath}
\usepackage{mathtools}
\mathtoolsset{mathic}
\makeatletter
\DeclareRobustCommand{\dotsc}{%
  \DN@{\ifx\@let@token;\@cdots\,%
       \else \ifx\@let@token.\@cdots\,%
       \else \extra@\@cdots \ifgtest@\,\fi
       \fi\fi}%
  \FN@\next@}
\makeatother
\usepackage{amssymb}
\usepackage{textgreek}
\usepackage{graphicx}
\newcommand\imp{\mathrel{\boldsymbol{\stixrightarrow}}}
\newcommand\equ{\mathrel{\boldsymbol{\stixrightleftarrows}}}
\newcommand\Imp{\mathrel{\rightarrow}}
\newcommand\IMP{\mathrel{\preccurlyeq}}
\newcommand\Equ{\mathrel{\leftrightarrow}}
\newcommand\rel{\mathrel{\vdash}}
\newcommand\homo{\mathrel{\rho}}
\newcommand\semicolon{\mathbin{\textsf{\textup{\textbf{;}}}}}
\newcommand\comma{\mathbin{\textsf{\textup{\textbf{,}}}}}
\newcommand\bigland{\mathop{\textstyle{\stixbigwedgeop}}\nolimits}
\newcommand\biglor{\mathop{\textstyle{\stixbigveeop}}\nolimits}
\newcommand\true{\mathord{\textstyle{\stixbigveeop}}}
\newcommand\false{\mathord{\textstyle{\stixbigwedgeop}}}
\newcommand\lt{\leqslant}
\newcommand\gt{\geqslant}
\renewcommand\exists{\biglor}
\renewcommand\forall{\bigland}
\renewcommand*\neg[1]{%
  \setbox0\hbox{$\mathaccent"0362{#1}^H$}%
  \setbox2\hbox{$\mathaccent"0362{\kern0pt#1}^H$}%
  \ifdim\ht0=\ht2 \overline{#1}\else \bar#1\fi
  }
\newcommand*\dng[1]{%
  \setbox0\hbox{$\mathaccent"0362{#1}^H$}%
  \setbox2\hbox{$\mathaccent"0362{\kern0pt#1}^H$}%
  \ifdim\ht0=\ht2 \overline{\overline{#1}}\else \bar{\bar#1}\fi
  }

\makeatletter
\newcommand\deutsch[1]{\leavevmode\vbox to0pt{{%
\vss\llap{\raise .65em%
\hbox{\normalfont\scriptsize\ttfamily#1}%
}}}}
\makeatother
\renewcommand\deutsch[1]{}
\newcommand\pv[2]{#2} % corrected error of the published version

\newenvironment{linedisplay}{\noindent\halign to\hsize\bgroup##\hfil\tabskip1em plus1fil&##\hfil&\hfil\llap{##}\tabskip0pt\cr}{\crcr\egroup}
\newcommand\Item[1]{\par\setbox0\hbox{#1 }\hangindent\wd0\noindent\box0\ignorespaces}

\newcounter{origpage}
\newcommand\bruch{\refstepcounter{origpage}\ensuremath{\mid^{\theorigpage}}}

\usepackage{paralist}
\makeatletter
\renewcommand\section{\@startsection{section}{1}{\parindent}{\medskipamount}{-\fontdimen2\font plus -\fontdimen3\font minus -\fontdimen4\font}{\normalfont\normalsize\fontseries{b}\selectfont}}
\makeatother
\usepackage{titlesec}
\titlelabel{\thetitle.\ }

\usepackage{amsthm}
\newtheoremstyle{jsl}{\medskipamount}{\medskipamount}{\slshape}{\parindent}{\scshape}{.}{ }{}
\theoremstyle{jsl}
\newtheorem{thm}{Theorem}

\makeatletter
\def\blfootnote{\xdef\@thefnmark{}\@footnotetext}
\makeatother

\title{Algebraic and logistic investigations on free lattices}
\author{Paul Lorenzen\textsuperscript{\textdied}\\\small translated by Stefan Neuwirth}
\date{}

\PassOptionsToPackage{urlcolor=magenta,citecolor=blue,anchorcolor=blue,linkcolor=blue,menucolor=cyan}{hyperref}%
\usepackage[pdftitle=,pdfauthor=,pdfdisplaydoctitle=true,pdflang=en-GB,bookmarksopen=false,breaklinks=true,unicode=true,plainpages=false,hyperindex=true,pdfstartview=FitH,colorlinks=true,pdfpagelabels=true,bookmarksnumbered=true]{hyperref}

\ifx\texorpdfstring\undefined\newcommand\texorpdfstring[2]{#1}\fi

\begin{document}
\setlength\abovedisplayskip{4pt plus 1pt minus 2pt}
\setlength\belowdisplayskip{\abovedisplayskip}
\maketitle

\blfootnote{Translation of “Algebraische und logistische Untersuchungen über freie Verbände”, The Journal of Symbolic Logic, 16(2), 81--106, published with the kind permission of
  Jutta Reinhardt, Paul Lorenzen’s daughter. Translator's note: in order to enhance the readability of the text for a public of lattice theorists as well as of logicians, we have replaced the signs~$<$, \textemdash, $\subset$, and~$\prec$ with the signs~$\lt$, $\rel$, $\subseteq$, and~$\IMP$, respectively; thirty misprints have been tacitly corrected (they are indicated in the \LaTeX\ code).}
It\label{intro-lattices} is well known that lattice theory was founded by Dedekind by means of\deutsch{an Hand} his ideal-theoretic investigations. It has turned out lately that the essential property of Dedekind's system of ideals lies in the fact that ideals form a semilattice (\S1). Ideal theory leads in this way\deutsch{so} to the question about all semilattices over an arbitrary preordered set~$M$. The simple answer is contained in~\S2. The question about all distributive lattices over~$M$ may be answered just\deutsch{ebenso} as simply. In both cases, among the semilattices vs.\ distributive lattices over~$M$, one is distinguished by the fact that all others are homomorphic to it. We call this distinguished semilattice vs.\ distributive lattice the free semilattice vs.\ free distributive lattice over~$M$.

One comes to new questions at investigating more special semilattices over~$M$. It is known e.g.\ that $M$~can always be extended into a complete boolean lattice. Is\label{q-complete-boolean} there, among these extensions, always the free complete boolean lattice over~$M$, which is distinguished by the fact that all others are homomorphic to it?

In~\S3, first the existence of the free orthocomplemented semilattice over~$M$ is proved. The method used here may also be followed to lead the proof of existence for the free countably complete boolean lattice over~$M$.

The significance\deutsch{Bedeutung} of the proofs of existence is not exhausted in pure lattice theory, but finds an important application in logistics. It is well known that the formalisation of logic has been---beside ideal theory---a further impulse for the development of lattice theory. Nevertheless\deutsch{trotzdem}, logistics were only able to exploit a modicum of lattice-theoretic results.

In~\S4 we however\deutsch{jedoch} show on a simple calculus of propositions how the question of freedom from contradiction and the decision problem is answered immediately by the proof of existence for free orthocomplemented semilattices.
\medskip

In part~II (\S\S5--8), the freedom from contradiction of ramified type logic including the axiom of infinity is being proved by the method of the proof of existence for the free count\-ably com\-plete boolean lattice. Knowledge of part~I (\S\S1--4) is however\deutsch{jedoch} not assumed.

By the fact that the basic thought of lattice theory is being used only implicitly, the proof of freedom from contradiction appears somehow as the continuation of the original approach by which Gentzen proved the freedom from contradiction of arithmetic without complete induction in his Ph.D.\ thesis. Freedom from contradiction results in fact as an immediate conclusion therefrom, that each theorem of the calculus may be deduced ``without detour''. The proof described here goes however\deutsch{jedoch} beyond Gentzen's proof, as\bruch{} the calculus whose freedom from contradiction is proved contains arithmetic including complete induction as part. This calculus is equivalent to the one used by Russell and Whitehead in the \emph{Principia mathematica} if the axiom of reducibility is removed there. As this axiom is not comprised, our calculus does not contain classical analysis, although the analytic modes of inference may still be represented in this calculus---with the restrictions required by ramified type theory.

The extension of Gentzen's approach to a so much richer calculus succeeds without addition of new means. Only the concept of deducibility without detour is extended by allowing certain induction rules in which a conclusion is inferred from infinitely many premisses. 

The progress with regard to the work of Fitch\footnote{F. B. Fitch, \emph{The consistency of the ramified} Principia, J. Symbolic Logic, vol.~3 (1938), pp.~140--149, and \emph{The hypothesis that infinite classes are similar}, ibid.,\ vol.~4 (1939), pp.~159--162.} lies in the constructive character of all inferences used. Only hereby does our proof fulfil the demands that have been addressed since Hilbert to a proof of freedom from contradiction.

In~\S5 the calculus whose freedom from contradiction is to be proved is presented. It will be called shortly the deductive calculus. It will be confronted in~\S6 to an inductive calculus that may be thought of as a specification of the concept of deducibility without detour. The inductive calculus is free from contradiction in a trivial way, so that for the freedom from contradiction of the deductive calculus one has to show that the inductive calculus is stronger than the deductive. This proof uses only inductions on formulae vs.\ theorems as auxiliary means, i.e.\ the fact that the concept of formula vs.\ theorem is defined constructively. In contrast\deutsch{Dagegen}, the so-called transfinite induction is not used.

By a little modification of the proof, it is established over and above in~\S8 that the axiom of reducibility is independent from the remaining axioms of the deductive calculus. In fact, the deductive calculus remains free from contradiction if countability of all sets is requested in addition. Cantor's diagonal procedure yields then in the extended calculus the refutability of the axiom of reducibility.

\section{Basic concepts.}
We gather first the basic concepts of the theory of semilattices.

Let $M$~be a set and $\lt$ a binary relation in~$M$. Let $a,b,\dotsc$ be the \pv{Element}{elements} of~$M$.

(A) $M$~is called ``preordered'' (w.r.t.~$\lt$) if holds:
\begin{linedisplay}
  1.&$a\lt a$.&\cr
  2.&${a\lt b\comma b\lt c}\quad\imp\quad{a\lt c}$.\footnote{I.e.\ $a\lt b$ and $b\lt c$ implies $a\lt c$.}&
\end{linedisplay}
\noindent Instead of $b\lt a$ we also write $a\gt b$. If $a\lt b$ and $a\gt b$ hold, then we write $a\equiv b$.\footnote{We do not assume ${a\equiv b}\imp {a=b}$. The concepts under (E)~and~(F) must hereby be defined somewhat differently than usual.} \ $\equiv$~is an equivalence relation.\bruch

(B) $M$~is called a ``semilattice'' (w.r.t.~$\lt$) if $M$~is preordered (w.r.t.~$\lt$) and if for each~$a$,~$b$ there is a~$c$ with:
\begin{linedisplay}
  3.1&$c\lt a$.&\cr
  3.2&$c\lt b$.&\cr
  3.3&${x\lt a}\comma{x\lt b}\quad\imp\quad{x\lt c}$.&
\end{linedisplay}
\noindent $c$~is uniquely determined (w.r.t.~$\equiv$). We write $c\equiv a\land b$. 

(C) $M$~is called a ``lattice'' (w.r.t.~$\lt$) if $M$~is a semilattice (w.r.t.~$\lt$) and simultaneously a semilattice (w.r.t.~$\gt$). If $c$~fulfils the conditions~3.1--3.3 with~$\gt$ instead of~$\lt$, then we write $c\equiv a\lor b$.

(D) $M$~is called a ``distributive lattice'' (w.r.t.~$\lt$) if $M$~is a lattice (w.r.t.~$\lt$) and if holds\\
4.\hfil ${a\land c\lt b}\comma{a\lt b\lor c}\quad\imp\quad{a\lt b}$.

(E) If $M$ vs.\ $M'$ is a preordered set (w.r.t.~$\lt$) vs.\ (w.r.t.~$\lt'$), then $M'$~is called a ``part'' of~$M$ if $M\pv{}{'}$~is a subset of~$M$ and if for each $a',b'\in M'$ holds $a'\lt'b'\equ a'\lt b'$.\footnote{I.e.\ $a'\lt'b'$ equivalent with $a'\lt b'$.} If $M$~is a semilattice vs.\ lattice, then $M$~is called a semilattice vs.\ lattice ``over~$M'$'' if $M'$~is a part of~$M$.

If $M$~is a semilattice vs.\ lattice over~$M'$, then $M$~is called a ``minimal'' semilattice vs.\ lattice over~$M'$ if $M$~does not contain a proper subset~$M_0$ for which holds:
\begin{linedisplay}
  (1)&$M'\subseteq M_0$.&\cr
  (2)&${a_0,b_0\in M_0}\comma{c\equiv a_0\land b_0}\quad\imp\quad{c\in M_0}$.&\cr
  (3)&${a_0,b_0\in M_0}\comma{c\equiv a_0\lor b_0}\quad\imp\quad{c\in M_0}$.&
\end{linedisplay}

(F) If $M$~and~$M'$ are preordered sets (w.r.t.~$\lt$), then a relation~$\homo$ between $M$~and~$M'$ is called a ``homomorphism'' from~$M$ into~$M'$ if holds:
\begin{linedisplay}
  [1]&To each~$a\in M$ there is an~$a'\in M'$ with $a\homo a'$.&\cr
  [2]&${a\homo a'_1}\comma{a'_1\equiv a'_2}\quad\imp\quad{a\homo a'_2}$.&\cr
  [3]&${a\homo a'}\comma{b\homo b'}\comma{a\lt b}\quad\imp\quad{a'\lt b'}$.&
\end{linedisplay}
\noindent If $M$~and~$M'$ are semilattices vs.\ lattices, then a ``homomorphism'' from~$M$ into~$M'$ is called a ``semilattice homomorphism'' vs.\ ``lattice homomorphism'' if holds:
\begin{linedisplay}
  [3.1]&${a\homo a'}\comma{b\homo b'}\quad\imp\quad{{a\land b}\homo{a'\land b'}}$.&\cr
  [3.2]&${a\homo a'}\comma{b\homo b'}\quad\imp\quad{{a\lor b}\homo {a'\lor b'}}$.&
\end{linedisplay}
\noindent(${a\land b}\homo {a'\land b'}$ means ${c\equiv a\land b}\comma {c'\equiv a'\land b'}\imp {c\homo c'}$.)

An homomorphism~$\homo$ from~$M$ into~$M'$ is called an ``isomorphism'' from~$M$ into~$M'$ if holds:\\\relax 
[4]\hfil ${a\homo a'}\comma {b\homo b}'\comma {a'\lt b'}\quad\imp\quad {a\lt b}$.

An homomorphism vs.\ isomorphism from~$M$ into~$M'$ is called an homomorphism vs.\ isomorphism from~$M$ ``onto~$M'$'' if holds:\\\relax
[5]\hfil To each~$a'\in M'$ there is an~$a\in M$ with $a\homo a'$.

$M'$ is called ``homomorphic'' vs.\ ``isomorphic'' to~$M$ if there is an homomorphism vs.\ isomorphism from~$M$ onto~$M'$. If $M'$~is homomorphic vs.\ isomorphic to~$M$, and if $M_0$~is a part of $M$~and~$M'$, then $M'$~is called homomorphic vs.\ isomorphic to~$M$ ``over~$M_0$'' if there is an homomorphism from~$M$ onto~$M'$ such that for each $a_0\in M_0$ holds $a_0\homo a_0$.\bruch

\section{Free semilattices and distributive lattices.}\label{sec:free}
Let $M$~be a preordered set. The minimal semilattices over~$M$ may be characterised by relations in~$M$ with a finite number of places, as show the following theorems.

\begin{thm}\label{thm1}
  If \(H\)~is a semilattice over~\(M\), then for the relation in~\(M\) defined by
  \[{a_1,\dotsc,a_n\rel b}\equ {a_1\land\dotsb\land a_n\lt b}\]
  holds:
\begin{linedisplay}
  $1$.&$a\rel a$.&\cr
  $2$.&${a_1,\dotsc,a_n\rel b}\imp {a_1,\dotsc,a_n,c\rel b}$.&\cr
  $3$.&${a_1,\dotsc,a_n\rel b}\imp {a_1,\dotsc,a_{i+1},a_i,\dotsc,a_n\rel b}$.&\cr
  $4$.&${a_1,\dotsc,a_n\rel c}\semicolon {a_1,\dotsc,a_n,c\rel b}\imp {a_1,\dotsc,a_n\rel b}$.&
\end{linedisplay}
\end{thm}

\begin{thm}
  If \(M\)~is preordered, then the relation defined by
  \[{a_1,\dotsc,a_n\rel b}\equ{(\text{there is an \(a_i\) with \(a_i\lt b\)})}\]
  fulfils the conditions \(1\).--\(4\).\ of theorem~\(1\).
\end{thm}

\begin{thm}
  To each relation \(a_1,\dotsc,a_n\rel b\) in~\(M\),\footnote{$M$~is preordered by the relation $a\rel b$.} that fulfils the conditions \(1\).--\(4\).\ of theorem~\(1\), there is an (up to isomorphy over~\(M\)) uniquely determined minimal semilattice over~\(M\) for which holds
  \[{a_1,\dotsc,a_n\rel b}\equ {a_1\land\dotsb\land a_n\pv{---}{\lt} b}.\]
\end{thm}

We call the semilattice associated to the relation of theorem~2 according to theorem~3 the ``free'' semilattice over~$M$.

\begin{thm}
  If \(H\)~is the free semilattice over~\(M\), then each minimal semilattice over~\(M\) is homomorphic\footnote{I.e.\ there is a semilattice homomorphism.} to~\(H\) over~\(M\).\label{thm4}
\end{thm}

The proofs of these theorems are so simple that we omit them.

For the proof of theorem~3, one forms the set~$H$ of all finite sequences $a_1\land\dotsb\land a_n$ out of elements of~$M$ and defines in~$H$ a preorder~$\lt$ by
\[{a_1\land\dotsb\land a_n\lt b_1\land\dotsb\land b_m}\equ{(\text{for each~$b_i$, \ $a_1,\dotsc, a_n\rel b_i$})}\text.\]
$H$~is the sought-after\deutsch{gesuchte} semilattice.

Minimal distributive lattices over~$M$ may be characterised just\deutsch{ebenso} as simply as semilattices.
\begin{thm}
  If \(V\)~is a distributive lattice over~\(M\), then for the relation in~\(M\) defined by
  \[{a_1,\dotsc, a_m\rel b_1,\dotsc, b_n}\equ {a_1\land\dotsb\land a_m\lt b_1\lor\dotsb\lor b_n}\]
  holds:
\begin{linedisplay}
  $1$.&$a\rel a$.&\cr
  $2$.&${a_1,\dotsc,a_m\rel b_1,\dotsc, b_n}\imp {a_1,\dotsc,a_m,c\rel b_1,\dotsc, b_n}$.&\cr
  &${a_1,\dotsc,a_m\rel b_1,\dotsc, b_n}\imp {a_1,\dotsc,a_m\rel c,b_1,\dotsc, b_n}$.\bruch&\cr
  $3$.&${a_1,\dotsc,a_m\rel b_1,\dotsc, b_n}\imp {a_1,\dotsc,a_{i+1},a_i,\dotsc,a_m\rel b_1,\dotsc, b_n}$.&\cr
  &${a_1,\dotsc,a_m\rel b_1,\dotsc, b_n}\imp {a_1,\dotsc,a_m\rel b_1,\dotsc,b_{i+1},b_i,\dotsc, b_n}$.&\cr
  $4$.&${a_1,\dotsc,a_m,c\rel b_1,\dotsc, b_n}\semicolon {a_1,\dotsc,a_m\rel c,b_1,\dotsc, b_n}$&\cr
  &&${}\imp {a_1,\dotsc,a_m\rel b_1,\dotsc, b_n}$.\cr
\end{linedisplay}
\end{thm}

\begin{thm}
If \(M\)~is preordered, then the relation defined by   
  \[{a_1,\dotsc,a_m\rel b_1,\dotsc, b_n}\equ{(\text{there is an~\(a_i\) and a~\(b_j\) with \(a_i\lt b_j\)})}\]
  fulfils the conditions \(1\).--\(4\).\ of theorem~\(5\).
\end{thm}

\begin{thm}\label{thm7}
  To each relation \(a_1,\dotsc,a_m\rel b_1,\dotsc, b_n\) in~\(M\),\footnotemark[5] that fulfils the conditions \(1\).--\(4\).\ of theorem~\(5\), there is an (up to isomorphy over~\(M\)) uniquely determined minimal distributive lattice over~\(M\) for which holds
  \[{a_1,\dotsc,a_m\rel b_1,\dotsc, b_n}\equ {a_1\land\dotsb\land a_m\lt b_1\lor\dotsb\lor b_n}\text.\]
\end{thm}

We call the distributive lattice associated to the relation of theorem~6 according to theorem~$7$ the ``free'' distributive lattice over~$M$.

\begin{thm}
  If \(V\)~is the free distributive lattice over~\(M\), then each minimal distributive lattice over~\(M\) is homomorphic\footnote{I.e.\ there is a lattice homomorphism.} to~\(V\) over~\(M\).
\end{thm}

We may again omit the proofs.

For the proof of theorem~7, one forms first---as at\deutsch{bei} proving theorem~3---the set~$H$ of finite sequences $a_1\land\dotsb\land a_n$ out of elements of~$M$. If $\alpha=a_1\land\dotsb\land a_m$ and $\beta_i=b_{i1}\land\dotsb\land b_{in_i}$ are elements of~$H$, then one sets $\alpha\rel\beta_1,\dotsc,\beta_n$ if $a_1,\dotsc,a_m\rel b_{1j_1},\dotsc, b_{nj_n}$ holds for each $j_1,\dotsc,j_n$. This relation fulfils conditions that correspond to~1.--4.\ of theorem~1. One forms then---correspondingly as at\deutsch{bei} proving theorem~3---the set~$V$ of finite sequences $\alpha_1\lor\dotsb\lor\alpha_n$ out of elements of~$H$ and defines a preorder~$\lt$ in~$H$ by
\[{\alpha_1\lor\dotsb\lor\alpha_m\lt\beta_1\lor\dotsb\lor \beta_n}\equ{(\text{for each~$\alpha_i$, \ $\alpha_i\rel\beta_1,\dotsc,\beta_n$})}\text.\]
$V$~is the sought-after\deutsch{gesuchte} lattice.\label{sec:free:end}

\section{Free orthocomplemented semilattices.}
After the general semilattices and the distributive lattices we investigate now more special semilattices.

A preordered set (w.r.t.~$\lt$) is called ``bounded'' if there are elements $0$~and~$1$ in~$M$ with $0\lt a\lt1$ for each~$a\in M$. $0$~vs.~1 is called ``zero element'' vs.\ ``unit element'' of~$M$. 

A bounded semilattice is called ``orthocomplemented'' if to each~$c$ there is a $d$ with
\[{a\land c\lt0}\equ {a\lt d}\text{ (for all~$a$).}\label{orthocomplementation}\]
We write $d\equiv\neg c$.

\fontdimen8\textfont3=.3pt
If $H$~is an orthocomplemented semilattice over~$M$, then $H$~is called a ``mini-\bruch{}mal''  orthocomplemented semilattice over~$M$ if $H$ contains no proper subset~$H'$ for which holds:
\begin{linedisplay}
  i.&$M\subseteq H'$&\cr
  ii.&${a'\in H'}\comma{b'\in H'}\comma{c\equiv a'\land b'}\quad\imp\quad {c\in H'}$.&\cr
  iii.&${c'\in H'}\comma {d\equiv\neg{c'}}\quad\imp\quad {d\in H'}$.&
\end{linedisplay}
If $H$~and~$H'$ are orthocomplemented semilattices, then a semilattice homomorphism~$\homo$ from~$H$ into~$H'$ is called an ``orthocomplemented'' semilattice homomorphism if holds
\[{a\homo a'}\imp{\neg a\homo\neg{a'}}\text.\]

Correspondingly to theorem~4 and theorem~8 we define now: an orthocomplemented semilattice over~$M$ is called a ``free'' orthocomplemented semilattice over~$M$ if each minimal orthocomplemented semilattice over~$M$ is homomorphic\footnote{I.e.\ there is an orthocomplemented semilattice homomorphism.} to~$H$. The free orthocomplemented semilattice over~$M$ is uniquely determined up to isomorphy over~$M$.

We prove below\deutsch{im folgenden} the existence of the free orthocomplemented semilattice over an arbitrary preordered set by construction. We prove more precisely:
\begin{thm}
  If \(M\)~is a bounded preordered set, and \(a_1,\dotsc,a_n\rel b\) a relation that fulfils the conditions
  \[{a\rel b}\equ {a\lt b}\leqno{1.}\]
  and \(2\).--\(4\). of theorem~\(1\), then there is an orthocomplemented semilattice~\(H\) over~\(M\) for which holds that for arbitrary elements \(a_1,\dotsc,a_n,b\) out of~\(M\),
  \[{a_1,\dotsc,a_n\rel b}\equ {a_1\land\dotsb\land a_n\lt b},\]
  and to which each minimal orthocomplemented semilattice that fulfils these conditions is homomorphic\footnotemark[8] over~\(M\).
\end{thm}

Let $H_0$~be the set~$M$. Let $H'_i$ be set of two-term sequences~$(a\land b)$ out of elements $a,b$ of~$H_i$, and let $H_{i+1}=H_i\cup\neg{H_i}\cup H'_i$. Let $H$ be the union of the sets~$H_i$ ($i=0,1,\dotsc$).

The elements of~$H$ we call shortly ``formulae,'' the elements of~$M$ ``prime formulae''. Then the following ``formula induction'' holds: if a claim holds
\begin{linedisplay}
  1.&for each prime formula,&\cr
  2.&for~$a\land b$, $\neg a$ if for~$a$,~$b$,&
\end{linedisplay}
\noindent then it holds for each formula.

In~$H$ we define constructively a relation~$\lt$ by:
\begin{compactenum}[(1)]
\item For prime formulae $a_1,\dotsc,a_n,b$, \ $a_1\land\dotsb\land a_n\lt b$ if $a_1,\dotsc,a_n\rel b$.
\item Structure rules. If $a_1\land\dotsb\land a_n\lt b$ holds, then a valid relation arises again by the following structure changes to the left formula: association,\bruch{}  i.e.\ the grouping by brackets may be changed;\footnote{We have hence omitted brackets beforehand.} transposition, i.e.\ the \pv{Reinhenfolge}{sequential arrangement} may be changed\pv{,}{.}
\end{compactenum}
\begin{linedisplay}
  (3.1)&${a\lt b}\imp {a\land c\lt b}$.&\cr
  (3.2)&${c\lt a}\comma {c\lt b}\imp {c\lt a\land b}$.&\cr
  (3.3)&${a\land b\lt0}\imp {a\lt\neg b}$.&\cr
  (3.4)&${a\lt b}\imp {a\land\neg b\lt c}$.&
\end{linedisplay}
\noindent In order to achieve that the unit element~$1$ of~$M$ also becomes the unit element of~$H$, we moreover\deutsch{ferner} set that these rules are also to hold if a formula~$1\land x$ is replaced by~$x$. Thus\deutsch{also} the rules (3.1), (3.3), and~(3.4) include:
\[\begin{gathered}
  {1\lt b}\imp {c\lt b}\text.\\
  {b\lt0}\imp {1\lt\neg b}\text.\\
  {1\lt b}\imp {\neg b\lt c}\text.
\end{gathered}\]
If $a\lt b$ holds, then we call the formula pair~$a,b$ a ``theorem.'' The theorems $a_1\land\dotsb\land a_n\lt b$ out of prime elements $a_1,\dotsc,a_n$, and~$b$ we call ``prime theorems.'' In the rules~(3) (and correspondingly in~(2)) we call the formula pairs to the left of~$\imp$ the ``premisses'', and the formula pair to the right of~$\imp$ the ``conclusion.'' The following ``theorem induction'' holds: if a claim holds
\begin{compactenum}[1.]
\item for each prime theorem,
\item for the conclusion of each rule whose premisses are theorems if for these pre\-mis\-ses,
\end{compactenum}
then it holds for each theorem.

We now show that $H$~is an orthocomplemented semilattice. For this we have to prove
\begin{linedisplay}
  (4.1)&$c\lt c$.&\cr
  (4.2)&$0\lt c\lt1$.&\cr
  (4.3)&${a\lt c}\comma {c\lt b}\imp {a\lt b}$.&\cr
  (5.1)&${a\lt b_1\land b_2}\imp {a\lt b_1}$.&\cr
  (5.2)&${a\lt b_1\land b_2}\imp {a\lt b_2}$.&\cr
  (5.3)&${a\lt\neg b}\imp {a\land b\lt0}$.&
\end{linedisplay}

We first prove~(4.1) by formula induction. For prime formulae~$c$ holds~$c\lt c$ because of~$c\rel c$. From $c_1\lt c_1$ and $c_2\lt c_2$ follows $c_1\land c_2\lt c_1$ and $c_1\land c_2\lt c_2$, thus $c_1\land c_2\lt c_1\land c_2$. From $c\lt c$ follows moreover\deutsch{ferner} $c\land\neg c\lt0$, and then $\neg c\lt\neg c$.

(4.2) follows as well\deutsch{ebenso} by formula induction. For prime formulae~$c$ holds~$0\lt c$ because of~$0\rel c$. From $0\lt c_1$ and $0\lt c_2$ follows $0\lt c_1\land c_2$. $0\lt\neg c$ holds because of ${0\lt0}\imp{0\land c\lt0}$, and $c\lt1$ follows directly from (3.1)~${1\lt1}\imp {c\lt1}$.

(4.3) we prove---because of the difficulty---last.

For the proof of~(5.1) and~(5.2) we use theorem induction. The induction claim states that for each theorem that has the form $a\lt b_1\land b_2$, also $a\lt b_1$ and $a\lt b_2$ are theorems. For prime theorems there is nothing to prove as there are no prime theorems of the form $a\lt b_1\land b_2$. Let now $a\lt b_1\land b_2$~be conclusion of a rule whose premisses are theorems, and let the claim hold for the premisses. There are then only the following possibilities:
\begin{compactenum}[(a)]
\item $a\lt b_1\land b_2$ is conclusion of a structure rule with a premiss~$a'\lt b_1\land b_2$ in which $a$ arises by a structure change out of~$a'$.\bruch
\item One has $a=a_1\land a_2$, and $a\lt b_1\land b_2$~is conclusion of rule~(3.1) with the premiss~$a_1\lt b_1\land b_2$.
\item $a\lt b_1\land b_2$~is conclusion of rule~(3.2) with the premisses~$a\lt b_1$ and~$a\lt b_2$.
\item One has $a=a_1\land\neg a_2$, and $a\lt b_1\land b_2$~is conclusion of rule~(3.4) with the premiss~$a_1\lt a_2$.
\end{compactenum}

In case~(a) holds according to the induction hypothesis $a'\lt b_1$ and~$a'\lt b_2$, from which~$a\lt b_1$ and~$a\lt b_2$ arise by a structure rule.

In case~(b) holds according to induction hypothesis $a_1\lt b_1$ and~$a_1\lt b_2$, from which according to~(3.1) follows $a_1\land a_2\lt b_1$ and $a_1\land a_2\lt b_2$.

In case~(c) holds $a\lt b_1$ and $a\lt b_2$, as the premisses are theorems.

In case~(d) follows from~$a_1\lt a_2$ according to~(3.4) $a_1\land\neg a_2\lt b_1$ and $a_1\land\neg a_2\lt b_2$.

Thereby (5.1) and (5.2) are proved by theorem induction. The single proof steps are all trivial, and may therefore be skipped in later similar cases.

Next (5.3) results from such a trivial theorem induction.

Of~(4.3) the special case $a\lt0\imp a\lt b$ may be proved immediately by formula induction on~$b$ and theorem induction on~$a$.

For the general case we need three lemmas:
\begin{linedisplay}
  (6)&$a\land\neg{\neg c\land d}\lt b\imp a\land c\lt b$.\footnote{By our agreement on the replacement of~$1\land x$ by~$x$, (6) includes: $a\pv{\lor}{\land}\dng{c}\lt b\imp a\land c\lt b$, $\dng{c}\lt b\imp c\lt b$. The rules arising by structure changes are---in order to abbreviate---not specified explicitly.}&
\end{linedisplay}
\noindent Proof by theorem induction.
\begin{linedisplay}
  (7)&If $a\land\neg c\lt p$ holds for a prime formula~$p$, then holds $a\lt p$ or $a\land \neg c\lt0$.&
\end{linedisplay}
\noindent Proof by theorem induction.
\begin{linedisplay}
  (8)&$a\land c\land c\lt b\imp a\land c\lt b$\label{contraction}.\footnote{These claims include: $a\land\neg c\lt c\imp a\lt c$, $\neg c\lt c\imp 1\lt c$, $c\land c\lt b\imp c\lt b$.}&
\end{linedisplay}
\noindent We use a formula induction on~$c$. For prime formulae~$c$, (8)~follows by theorem induction. If (8)~holds for $c_1$~and~$c_2$, then naturally also for~$c_1\land c_2$. For formulae~$\neg c$ we prove
\[a\land\neg c\land\neg c\lt b\imp a\land\neg c\lt b\]
by theorem induction. All steps are trivial, except in the case in which $a\land\neg c\land\neg c\lt b$ is conclusion of rule~(3.4) with the premiss~$a\land\neg c\lt c$. Let $c=p_1\land\dotsb\land p_m\land\neg c_1\land\dotsb\land\neg c_n$ with prime formulae~$p_\mu$. We have as induction hypothesis of our formula induction the validity of~(8) for each~$c_\nu$ ($\nu=1,\dotsc,n$). From $a\land\neg c\lt p_\mu$ follows according to~(7) $a\land\neg c\lt0$ or $a\lt p_\mu$ for each~$p_\mu$. From~$a\land\neg c\lt\neg c_\nu$ follows according to~(6) $a\land c_\nu\lt\neg c_\nu$; thus according to~(5.3) $a\land c_\nu\land c_\nu\lt0$, and therefore\deutsch{daher} $a\land c_\nu\lt0$, thus~$a\lt\neg c_\nu$. Together follows $a\land\neg c\lt0$ or $a\lt p_1\land\dotsb\land p_m\land\neg c_1\land\dotsb\land\neg c_n$; thus in each case $a\land\neg c\lt0$ and $a\land\neg c\lt b$.

By the aid of~(8) we can now instead of~(4.3) even prove
\begin{linedisplay}
  (9)&$a_1\lt c\comma a_2\land c\lt b\imp a_1\land a_2\lt b$&
\end{linedisplay}
\noindent by formula induction for each formula~$c$.

For prime formulae~$c$ proof by theorem induction.\bruch

If (9)~holds for $c_1$~and~$c_2$, then also for~$c_1\land c_2$, for from $a_1\lt c_1\land c_2$ and $a_2\land c_1\land c_2\lt b$ follows according to (5.1),~(5.2) $a_1\lt c_1$ and $a_1\lt c_2$; thus $a_1\land a_1\land a_2\lt b$, i.e.\ $a_1\land a_2\lt b$.

Now let (9)~hold for~$c$. We prove the validity for~$\neg c$ by theorem induction for all theorems~$a_2\land\neg c\lt b$. All steps are trivial, except in the case in which $a_2\land\neg c\lt b$~is conclusion of rule~(3.4) with the premiss~$a_2\lt c$. From~$a\lt\neg c$ follows according to~(5.3) and the already proved special case of~(4.3), $a_1\land c\lt b$; thus $a_1\land a_2\lt b$ according to induction hypothesis because of~$a_2\lt c$. 

Thereby\deutsch{damit} (9)~is proved and especially~(4.3).

For the proof of theorem~9 we have to note now that for elements $a_1,\dotsc,a_n,b$ out of~$M$, $a_1,\dotsc,a_n\rel b$ holds exactly if $a_1\land\dotsb\land a_n\lt b$; for out of prime theorems do only through\deutsch{durch} structure rules and rule~(3.1) arise again prime theorems. But these rules follow from conditions~1.--4.\ of theorem~1. Now only the verification that each minimal orthocomplemented semilattice~$H'$ over~$M$ for which
\[a_1,\dotsc,a_n\rel b\equ a_1\land\dotsb\land a_n\lt b\]
holds is homomorphic to~$H$ over~$M$ is lacking. We define for this inductively a relation~$\homo$ between $H$~and~$H'$ by:
\begin{linedisplay}
  (i)&$a\homo a$ for~$a\in M$.\label{inductiverelation}&\cr
  (ii)&$a\homo a'\comma b\homo b'\imp a\land b\homo a'\land b'$.&\cr
  (iii)&$a\homo a'\imp\neg a\homo\neg a'$.&
\end{linedisplay}
\noindent$\homo$ is an homomorphism, for holds
\[a\homo a'\comma b\homo b'\comma a\lt b\imp a'\lt b'\text,\]
as follows at once by theorem induction. The proof relies simply on the fact that rules~(2) and~(3) are valid for each orthocomplemented semilattice.

\section{Logistic application and complete lattices.}
The fact that the \label{calculuses}logic calculuses are semilattices or lattices permits a simple logistic application of free lattices.

A calculus whose formulae form an orthocomplemented semilattice arises in the following way.

We start with propositional variables and add all formulae $a,b,\dotsc$ that can be formed out of them by use of the conjunction sign~$\land$ and the negation sign~$\neg{\phantom{a}}$. Let $a\land b$~mean the proposition~``$a$~and~$b$,'' $\neg a$ the proposition~``not~$a$.'' To the formulae we add moreover\deutsch{ferner} $0$~and~$1$. Let $0$ mean the~``false,'' $1$ the~``true.'' In the set~$H$ of formulae we define a relation~$\lt$. Let $a\lt b$ mean ``$a$ implies~$b$.'' Let hold
\[\begin{gathered}
  a\lt a\text,\\
  0\lt a\lt1\text.
\end{gathered}\]
Let moreover\deutsch{ferner} hold each relation that may be derived from this on the basis\deutsch{aufgrund} of the following rules:
\[\begin{aligned}
  a\lt c\comma c\lt b&\imp a\lt b\text.\\
  c\lt a\comma c\lt b&\imp c\lt a\land b\text{.\bruch}\\
  c\lt a\land b&\imp c\lt a\text.\\
  c\lt a\land b&\imp c\lt b\text.\\
  a\land b\lt0&\imp a\lt\neg b\text.\\
  a\lt\neg b&\imp a\land b\lt0\text.
\end{aligned}\]

According to this definition, the relation~$\lt$ is a preorder, and~$H$ is an orthocomplemented semilattice with respect to~$\lt$. For the set~$M$ of propositional variables including~$0$ and~$1$, a preorder is defined through\deutsch{durch} the relations
\[\begin{gathered}
  a\lt a\text,\\
  0\lt a\lt1\text.
\end{gathered}\]
According to the result of~\S3, there exists the free orthocomplemented semilattice~$H_0$ over~$M$. $H$~is isomorphic to~$H_0$ over~$M$, for $H$~and~$H_0$ consist of the same formulae, and each relation~$a\lt b$ that holds for~$H$ holds also for~$H_0$, as well\deutsch{ebenso} as the other way round. Thus~$M$ is a part of~$H$. From this follows immediately the freedom from contradiction\footnote{The freedom from contradiction means that no proposition is simultaneously true and false, i.e.\ for no formula~$a$ does $1\lt a$ and $a\lt0$ hold simultaneously.} of the calculus, for from~$1\lt a$ and~$a\lt0$ would follow~$1\lt0$; but this relation does not hold in~$M$.

Apart from\deutsch{außer} the freedom from contradiction, a decision procedure for~$H$ follows from the construction
of~$H_0$. In fact, the validity of a relation~$a\lt b$ in~$H_0$ is decidable, as obviously for each theorem only finitely many premisses are possible and a chain of premisses always stops after finitely many terms (the maximal number of steps is easy to estimate).

Also the freedom from contradiction of formalised theories in which one does not work with this propositional calculus but with the classical calculus results by this means, if one considers instead of the orthocomplemented semilattices the count\-ably com\-plete boolean lattices.

A ``boolean lattice'' is a distributive lattice in which to each element~$c$ there is an element~$\neg c$ with $c\land\neg c\lt0$ and $c\lor\neg c\gt1$.

A semilattice~$H$ (with respect to~$\lt$) is called ``count\-ably com\-plete'' if for each countable subset~$N$ of~$H$ there is an element~$c$ of~$H$ such that holds:
\[\begin{gathered}
  a\in N\imp c\lt a\text.\\
  \text{(for each }a\in N\text{, }x\lt a\text)\imp x\lt c\text.
\end{gathered}\]
We then write $c=\bigland N$.

A lattice~$V$ is called ``count\-ably com\-plete'' if $V$~is a count\-ably com\-plete semilattice with respect to~$\lt$ and w.r.t.~$\gt$.

Classical number theory e.g., which for each formula $a(x)$ in which occurs\deutsch{vorkommen} a free individual variable~$x$ also contains the formulae $\forall_xa(x)$ and $\exists_xa(x)$ ($\forall_xa(x)$ means ``for each~$x$, $a(x)$,'' $\exists_xa(x)$ means ``for at least one~$x$, $a(x)$'') is indeed not a  count\-ably com\-plete lattice; it does not contain for each countable subset~$N$ of formulae e.g.\ the conjunction $\forall N$, but only for the sets $N=\{a(1),a(2),\dotsc\}$, and then one has $\forall N\equiv\forall_xa(x)$.\bruch

Nevertheless\deutsch{trotzdem}, the proof of existence for the free count\-ably com\-plete boolean lattice over any preordered set and the proof of freedom from contradiction for the classical calculus are so alike that, to avoid repetitions, we only sketch here the proof of existence. All details may be extracted from the proof of freedom from contradiction undertaken in part~II (\S\S5--8).

\begin{thm}
  If \(M\)~is a bounded preordered set and \(a_1,\dotsc,a_{\pv{n}{m}}\rel b_1,\dotsc,b_n\) a relation that satisfies the conditions
\begin{linedisplay}
  $1$.&$a\rel b\equ a\lt b$&
\end{linedisplay}
\noindent and \(2\).--\(4\).\ of theorem~\(5\), then there is a count\-ably com\-plete boolean lattice~\(V\) over~\(M\) for which holds that for elements \(a_1,\dots,a_m,b_1,\dots,b_n\in M\)
\[a_1,\dotsc,a_m\rel b_1,\dotsc,b_n\equ a_1\land\dotsb\land a_m\lt b_1\lor\dotsb\lor b_n\text,\]
and that each minimal\footnote{A count\-ably com\-plete boolean lattice~$V$ is called minimal over~$M$ if $V$~does not contain a proper subset~$V_0$ for which holds:
  \begin{compactenum}[1.]
  \item $M\subseteq V_0$.
  \item $a,b\in V_0\comma c\equiv a\land b\imp c\in V_0$.
  \item $a,b\in V_0\comma c\equiv a\lor b\imp c\in V_0$.
  \item $a\in V_0\comma c\equiv\neg a\imp c\in V_0$.
  \item $N\subseteq V_0\comma N\text{ countable}\comma c\equiv\forall N\imp c\in V_0$.
  \item $N\subseteq V_0\comma N\text{ countable}\comma c\equiv\exists N\imp c\in V_0$.%
  \end{compactenum}} count\-ably com\-plete boolean lattice~$V'$ over~$M$ that fulfils these conditions is homomorphic\footnote{I.e.\ there is a lattice homomorphism~$\homo$ for which holds:
  \begin{compactenum}[1.]
  \item $a\homo a'\imp\neg a\homo\neg a'$.
  \item If, for a countable subset~$N$, $\homo$~is an homomorphism from~$N$ onto~$N'$, then holds $\forall N\homo\forall N'$ and $\exists N\homo\exists N'$.%
  \end{compactenum}} to~$V$ over~$M$.
\end{thm}
 
First a set~$V$ will be defined constructively, for which holds:
\[\begin{gathered}
  M\subseteq V\text.\\
  a,b\in V\imp a\land b\in V\text.\\
  a,b\in V\imp a\lor b\in V\text.\\
  a\in V\imp\neg a\in V\text.\\
  N\subseteq V\comma N\text{ countable}\imp\forall N\in V\text.\\
  N\subseteq V\comma N\text{ countable}\imp\exists N\in V\text.
\end{gathered}\]
In doing so, let differently designated elements always be different.

In~$V$ a relation is defined constructively by:
\begin{linedisplay}
  [1]&$a_1,\dotsc,a_m\rel b_1,\dotsc,b_n\imp a_1\land\dotsb\land a_m\lt b_1\lor\dotsb\lor b_n$.&\cr
  \noalign{\Item{[2]}For the left and right formula of $a_1\land\dotsb\land a_m\lt b_1\lor\dotsb\lor b_n$ will be allowed as structure change, apart from\deutsch{außer} association and transposition, also contraction, i.e.\ of two equal elements one may be omitted.}
  [3]&$1\land x$ and $0\lor x$ may always be replaced by~$x$.\bruch&\cr
  [3.1]&$a\lt b\imp a\lt b\lor c$.&\cr
  [3.2]&$a\lt b\imp a\land c\lt b$.&\cr
  [3.3]&$a_1\land c\lt b\comma a_2\land c\lt b\imp(a_1\lor a_2)\land c\lt b$.&\cr
  [3.4]&$a\lt b_1\lor c\comma a\lt b_2\lor c\imp a\lt(b_1\land b_2)\lor c$.&\cr
  [3.5]&$a\lt b\lor c\imp a\land\neg c\lt b$.&\cr
  [3.6]&$a\land c\lt b\imp a\lt b\lor\neg c$.&
\end{linedisplay}
\noindent For~$N\subseteq V$, $N$~countable:
\begin{linedisplay}
  [3.7]&$c\in N\comma a\land c\lt b\imp a\land\forall N\lt b$.&\cr
  [3.8]&$c\in N\comma a\lt b\lor c\imp a\lt b\lor\exists N$.&\cr
  [3.9]&$(\text{For each~$x\in N$, $a\land x\lt b$})\imp a\land\exists N\lt b$.&\cr
  [3.10]&$(\text{For each~$x\in N$, $a\lt b\lor x$})\imp a\lt b\lor\forall N$.\label{booleanomegarule}&
\end{linedisplay}
\noindent By the aid of formula and theorem inductions one has then to prove that $V$~is w.r.t.~$\lt$ a count\-ably com\-plete boolean lattice.

Instead of the relations~(5.1),~(5.2) to be proved steps in\deutsch{treten} now:
\vskip\abovedisplayskip{\offinterlineskip\halign{\hskip\parindent\hfil\strut$#$.&\quad\vrule width1pt depth\abovedisplayskip#\quad&\hfil$#$\cr
    a\lt(b_1\land b_2)\lor c\imp a\lt b_1\lor c&&x\in N\comma a\lt b\lor\forall N\imp a\lt b\lor x\text.\cr
    a\lt(b_1\land b_2)\lor c\imp a\lt b_2\lor c&&\cr
    (a_1\lor a_2)\land c\lt b\imp a_1\land c\lt b&&x\in N\comma a\land\exists N\lt b\imp a\land x\lt b\text.\cr
    (a_1\lor a_2)\land c\lt b\imp a_2\land c\lt b&depth2pt&\cr
}}\vskip\belowdisplayskip
\noindent Instead of~(5.3) steps in\deutsch{treten}:
\[\begin{gathered}
  a\lt\neg b\lor c\imp a\land b\lt c\text.\\
  a\land\neg b\lt c\imp a\lt b\lor c\text.
\end{gathered}\]
Instead of~(4.3) vs.\ (9) steps in\deutsch{treten}
\[a_1\lt c\lor b_1\comma a_2\land c\lt b_2\quad\imp\quad a_1\land a_2\lt b_1\lor b_2\text.\label{boolecutrule}\]
The special case $a\lt0\imp a\lt b$ results immediately, $1\lt b\imp a\lt b$ as well\deutsch{ebenso}.
(8) (and therefore\deutsch{daher} also~(6) and~(7)) is dispensable here.

The essential difference with regard to the construction in theorem~9 lies in the fact that in~$V$, the validity of a relation~$a\lt b$ cannot be decided in general, for in the rules~[3.9] and~[3.10] occur\deutsch{vorkommen} infinitely many premisses.

The structure rule of contraction, that was provable in~\S3, must be assumed here, as shows the following example. Let\label{counterexample} $M=\{1,\frac12,\dotsc,1/n,\dotsc,0\}$, and let $\lt$~be the order according to magnitude. $V$~contains for $N=\{1,\frac12,\dotsc,1/n,\dotsc\}$ the element~$a=\neg{\forall N}$; moreover\deutsch{ferner} for $N'=N\cup\{a\}$ the element~$\forall N'$. For each~$c\in N$ holds then $c\lt c$, from which follows
\[\begin{gathered}
  \forall N'\lt c\text,\qquad\forall N'\lt\forall N\text,\\
  \forall N'\land a\lt0\text,\\
  \forall N'\land\forall N'\lt0\text{.\bruch}
\end{gathered}\]
In contrast\deutsch{Dagegen}---without assuming contraction---$\forall N'\lt0$ does not hold, for $c\lt0$ holds for no~$c\in N$; and also
\[\neg{\forall N}\lt0\]
does not hold, as $1\lt\forall N$ does not hold. The example shows further that contraction also does not become provable if the rules~\pv{(3.7)}{[3.7]} and~\pv{(3.8)}{[3.8]} are replaced by
\[\begin{gathered}
  c_1,\dotsc,c_n\in N\comma a\land c_1\land\dotsb\land c_n\lt b\qquad\imp\qquad a\land\forall N\lt b\text,\\
  c_1,\dotsc,c_n\in N\comma a\lt b\lor c_1\lor\dotsb\lor c_n\qquad\imp\qquad a\lt b\lor\exists N\text.
\end{gathered}\]

By the hypothesis of contraction, the proof of theorem~10 simplifies considerably with regard to the proof of theorem~9. The proof of freedom from contradiction undertaken in part~II (\S\S5--8) requires in contrast\deutsch{dagegen} again additional considerations because of the use of free variables in the logic calculus.

From theorem~10 the freedom from contradiction of ramified type logic (incl.\ arithmetic) may be derived immediately in the following way:

We define a calculus~$Z_0$.\\
Numbers: $\alpha,\beta,\dotsc$.
\begin{linedisplay}
  (N1)&1.&\cr
  (N2)&$\alpha\imp\alpha+1$.&
\end{linedisplay}
\noindent Formulae: $a,b,\dotsc$.
\begin{linedisplay}
  (F1)&$\true$. (Interpretation: the true.)&\cr
  &$\false$. (Interpretation: the false.)&\cr
  (F2)&$\alpha>\beta$.&
\end{linedisplay}
\noindent Expressions: $A,B,\dotsc$.
\begin{linedisplay}
  &$a\lt b$ (if~$a$, then~$b$).&
\end{linedisplay}
\noindent Theorems:
\begin{linedisplay}
  (T1)&$a\lt b$ if interpretation true.&\cr
  (T2)&$a\lt a$.&\cr
  &$\false\lt a\lt\true$.&\cr
  (T3)&$a\lt c\comma c\lt b\imp a\lt b$.&
\end{linedisplay}
\noindent$Z_0$~is obviously free from contradiction in the sense that $\true\lt\false$ is not a theorem.

It is to be shown that this freedom from contradiction is conserved if free vs.\ bound number variables $x,y,\dotsc$ are added, the definition of formulae is extended by
\begin{linedisplay}
  (F3)&$x> y,\qquad x>\alpha,\qquad\alpha> x$,&\cr
  (F4)&if $a,b$, then $a\land b,a\lor b$ (and, or),&\cr
  (F5)&if $c$, then $\neg c$ (not),&\cr
  (F6)&if $a(x)$, then $\forall_xa(x)$, $\exists_xa(x)$ (for all, for some),&
\end{linedisplay}
\noindent and the definition of theorems by
\begin{linedisplay}
  (T4)&$a\lt b_1\comma a\lt b_2\equ a\lt b_1\land b_2$,&\cr
  (T5)&$a_1\lt b\comma a_2\lt b\equ a_1\lor a_2\lt b$,&\cr
  (T6)&$a\land c\lt b\imp a\lt b\lor\neg c$,&\cr
  (T7)&$a\lt b\lor c\imp a\land\neg c\lt b$,&\cr
  (T8)&$a\lt b(x)\equ a\lt\forall_xb(x)$ ($x$ not in $a$),&\cr
  (T9)&$a(x)\lt b\equ\exists_xa(x)\lt b$ ($x$ not in $b$),&\cr
  (T10)&$\text{(for each $\alpha$, $A(\alpha)$)}\equ A(x)$.&
\end{linedisplay}
\noindent We call the arising calculus~$K_0$.\bruch

The expressions of~$K_0$ without free variables form a partial calculus~$K_0'$ that contains~$Z_0$. In the definition of theorems of~$K_0'$ step in\deutsch{treten}, instead of~(T8)--(T10):
\[\begin{gathered}
  \text{(for each $\gamma$, $a\lt b(\gamma)$)}\equ a\lt\forall_xb(x)\text,\\
  \text{(for each $\gamma$, $a(\gamma)\lt b$)}\equ\exists_xa(x)\lt b\text.
\end{gathered}\]
Therefore $K_0'$~is contained isomorphically in the free count\-ably com\-plete boolean lattice over~$Z_0$, i.e.\ $Z_0$~is a partial calculus of~$K_0'$---in particular $K_0'$~and~$K_0$ are thus free from contradiction.

The freedom from contradiction follows as well\deutsch{ebenso} if in~$Z_0$ one admits as formulae apart from\deutsch{außer}~$\alpha>\beta$ in addition\deutsch{noch} $\alpha+\beta=\gamma$, $\alpha\cdot\beta=\gamma$ and the like.

In order to obtain the freedom from contradiction of ramified type logic (incl.\ arithmetic), we extend~$K_0$ by adding for each formula~$a(x)$ out of~$K_0$ the ``set''~$A=\hat xa(x)$, by extending the definition of formulae by
\begin{linedisplay}
  (F2*)&$x\in A,\qquad \alpha\in A$,&
\end{linedisplay}
\noindent and the definition of theorems by
\begin{linedisplay}
  (T1*)&$a(x)\lt x\in A,\qquad x\in A\lt a(x)$,&
\end{linedisplay}
\noindent and (T2), (T3), and (T10) also for the formulae (F2*). The calculus~$Z_1$ arising in this way\deutsch{so} is free from contradiction as $K_0$~is free from contradiction.

This freedom from contradiction is conserved by reason of\deutsch{aufgrund} the existence of the free count\-ably com\-plete boolean lattice over~$Z_1$ if free and bound set variables $X,Y,\dotsc$ are added, the definition of formulae is extended by
\begin{linedisplay}
  (F3*)&$x\in X,\qquad\alpha\in X$,&\cr
  (F6*)&if $a(X)$,\quad then $\forall_Xa(X)$,\quad $\exists_Xa(X)$,&
\end{linedisplay}
\noindent and (F4)--(F6) also for the formulae (F2*), (F3*), (F6*), and the definition of theorems by
\begin{linedisplay}
  (T8*)&$a\lt b(X)\equ a\lt\forall_Xb(X)$&($X$ not in~$a$),\cr
  (T9*)&$a(X)\lt b\equ\exists_Xa(X)\lt b$&($X$ not in $b$),\cr
  (T10*)&$(\text{for each $A$, $\mathfrak A(A)$)}\equ\mathfrak A(X)$.&
\end{linedisplay}
\noindent Iteration of this extension procedure yields the sought-after\deutsch{gesuchte} freedom from contradiction of ramified type logic incl.\ arithmetic.

\section{The deductive calculus of ramified type logic.}
Let $\lambda $~be a constructible ordinal number (e.g.\ $\omega$). By ``ordinal number'' we are always understanding below\deutsch{im folgenden} only the ordinal numbers~$\nu<\lambda$. We assume the knowledge of the signs for these ordinal numbers $0,1,2,\dotsc$, and are building up our calculus out of them by adding finitely many further individual signs:
\[\mathord(,\mathord),\mathord{\land},\mathord{\lor},\neg{\phantom{a}},\forall,\exists\text.\]

By ``sign'' we understand not only the individual signs, but also their compositions.
For communication we use for signs mostly $a,b,\dotsc$. $a=b$~means that~$a$ and~$b$ are signs of the same form. If $a$~is an individual sign, then $c(a)$~means a sign in which $a$~occurs\deutsch{vorkommen}. $c(b)$~means then the sign arising by substitution of~$a$ by~$b$. (At this\deutsch{hierbei}, each occurrence\deutsch{vorkommen} of~$a$ is to be substituted.) Also if $a$~is composite, we use this notation in the cases in which no misunderstanding is possible. We communicate sign pairs~$a,b$\bruch{} through\deutsch{durch} the letters $A,B,\dotsc$. Let $C(a)$~mean a sign pair in which $a$~occurs\deutsch{vorkommen} in at least one sign. Let $C(b)$~mean the sign pair arising by substitution of~$a$ by~$b$.

For the set-up of the deductive calculus we first define which signs we want to call ``types''. As ``type of $0$th~order'' we only take: $0$. As ``types of $\nu$th~order'' ($\nu>0$) we take: $\nu\,(\tau_1\tau_2\dotsm\tau_n)$ if $\tau_1,\dotsc,\tau_n$ are types of orders $\mu_1,\dotsc,\mu_n$ with $\mu_i<\nu$ ($i=1,\dotsc,n$).

For each type~$\tau$ we then form ``free variables'' $(\tau)_0, (\tau)_1,\dotsc$, and ``bound variables'' $((\tau))_0,\allowbreak((\tau))_1,\dotsc$. As signs of communication for free vs.\ bound variables we use $p,q,\dotsc$ vs.\ $x,y,\dotsc$. These letters we use possibly with indices. $p^\tau,q^\tau,\dotsc$ vs.\ $x^\tau,y^\tau,\dotsc$ mean always variables of type~$\tau$. Variables of type~$0$ we also call ``individual variables'', the variables of types of higher order also ``relation variables.''

Next we define which signs are to be called ``formulae''. As formula of $\nu$th order ($\nu>0$) we take:
\begin{linedisplay}
  (F1)&$0,1$.&\cr
  (F2)&$(p^\tau p_1^{\tau_1}\dotsm p_n^{\tau_n})$ if $\tau=\mu(\tau_1\dotsm\tau_n)$ and $\mu\leqq\nu$.&
\end{linedisplay}
\noindent These formulae are called the ``prime formulae.''

\Item{(F3)}With $a$, $b$, $c(p^\tau)$ also $(a\land b)$, $(a\lor b)$, $\neg a$, $\forall_{x^\tau}c(x^\tau)$, $\exists_{x^\tau}c(x^\tau)$ if $x^\tau$~does not occur\deutsch{vorkommen} in~$c(p^\tau)$ and for the order~$\mu$ of~$\tau$ holds $\mu<\nu$.

\noindent For the logical interpretation, $p^\tau p_1^{\tau_1}\dotsm p_n^{\tau_n}$ is to be read as ``the relation~$p^\tau$ is fulfilled by $p_1^{\tau_1},\dotsc,p_n^{\tau_n}$.'' $\mathord{\land}$, $\mathord{\lor}$, $\neg{\phantom{a}}$, $\forall$, $\exists$ is to be read as ``et,'' ``vel,'' ``non,'' ``omnes,'' ``existit.''

At\deutsch{bei} communicating formulae we omit the brackets as soon as this is possible without misunderstanding. For the communication of formulae~$(\neg a\lor b)$ vs.\ $(\neg a\lor b)\land(\neg b\lor a)$ we also use $a\Imp b$ vs.\ $a\Equ b$. We define last the concept of theorem for our calculus. In an only formal opposition to the classical calculus we are not distinguishing certain formulae as theorems, but formula pairs. We write $a\IMP b$ for communicating that the formula pair~$a,b$ is a theorem. For the logical interpretation, this is to be read as ``the proposition~$a$ implies the proposition~$b$.'' $0$~is the ``false,'' $1$~is the ``true.''

The concept of theorem is defined by:
\begin{linedisplay}
  (1a)&$c\IMP c,\qquad0\IMP c,\qquad c\IMP1$.&
\end{linedisplay}
\noindent(1b)\label{(1b)} For $\tau=\nu\,(\tau_1\dotsm\tau_n)$ and formulae of $\nu$th order $c(p_1^{\tau_1},\dotsc,p_n^{\tau_n})$,
\[1\IMP\exists_{x^\tau}\forall_{x_1^{\tau_1}}\dotsm\forall_{x_n^{\tau_n}}(x^\tau x_1^{\tau_1}\dotsm x_n^{\tau_n}\Equ c(x_1^{\tau_1},\dotsc,x_n^{\tau_n}))\text.\]
(1c) For $\tau=\nu\,(\tau_1\dotsm\tau_n)$ and formulae~$c(p^\tau)$,
\[\forall_{x_1^{\tau_1}}\dotsm\forall_{x_n^{\tau_n}}(p^\tau x_1^{\tau_1}\dotsm x_n^{\tau_n}\Equ q^\tau x_1^{\tau_1}\dotsc x_n^{\tau_n})\IMP c(p^\tau)\Imp c(q^\tau)\text.\]
(1d) For $\tau=1(00)$,
\[\makebox[\textwidth]{\(1\IMP\exists_{\!\!x^\tau}(\forall_{x^0}\neg{x^\tau x^0x^0}\land\forall_{x^0}\!\!\exists_{\!\!y^0}x^\tau x^0y^0\land\forall_{x^0}\!\forall_{y^0}\!\forall_{z^0}((x^\tau x^0y^0\land x^\tau y^0z^0)\Imp x^\tau x^0z^0))\).}\]
These theorems are called ``axioms.''

\noindent\halign to\hsize{#\hfil\tabskip1em plus1fil&\hfil$#$\tabskip3pt&$\imp$#&\tabskip1em plus1fil$#$\hfil\cr
(2a)&a\IMP c\comma c\IMP b&&a\IMP b\text{.\bruch}\cr%}\halign to\hsize{#\hfil\tabskip1em plus1fil&\hfil$#$\tabskip3pt&$\imp$#&\tabskip1em plus1fil$#$.\hfil\cr
(2b)&a\IMP b&&a\land c\IMP b\text.\cr
&a\IMP b&&c\land a\IMP b\text.\cr
&a\IMP b&&a\IMP b\lor c\text.\cr
&a\IMP b&&a\IMP c\lor b\text.\cr
&a\IMP b_1\comma a\IMP b_2&&a\IMP b_1\land b_2\text.\cr
&a_1\IMP b\comma a_2\IMP b&&a_1\lor a_2 \IMP b\text.\cr
(2c)&a\land c\IMP b&&a\IMP b\lor\neg c\text.\cr
&a\IMP b\lor c&&a\land\neg c\IMP b\text.\cr}
\noindent(2d) For formulae~$c(p^\tau)$ in which $x^\tau$ does not occur\deutsch{vorkommen},
\[\begin{aligned}
  a\IMP\forall_{x^\tau}c(x^\tau)&\imp a\IMP c(p^\tau)\text,\\
  \exists_{x^\tau}c(x^\tau)\IMP b&\imp c(p^\tau)\IMP b\text;
\end{aligned}\]
and, if $p^\tau$~does not occur\deutsch{vorkommen} in~$a$ nor~$b$,
\[\begin{aligned}
  a\IMP c(p^\tau)&\imp a\IMP\forall_{x^\tau}c(x^\tau)\text,\\
  c(p^\tau)\IMP b&\imp\exists_{x^\tau}c(x^\tau)\IMP b\text.
\end{aligned}\]
We call (2a)--(2d) the ``rules'' of the calculus. The propositions to the left of~$\imp$ are called the ``premisses,'' the proposition to the right of~$\imp$ is called the ``conclusion.''

The calculus is described by the complete definition of the concept of formula and theorem. We call it shortly the ``deductive calculus.''

Comments to the deductive calculus. (A)~If we restrict ourselves to ordinal numbers~$<\omega$, then our calculus is equivalent to the ``classical'' calculus of the \emph{Principia mathematica} if the axiom of reducibility is omitted there. A formula is classically deducible exactly if $1\IMP a$~holds in our calculus, and $a\IMP b$~holds exactly if the formula~$a\Imp b$ is classically deducible. We shall not undertake the proof of equivalence here, because it results easily from the equivalence of Gentzen's sequent calculus with the classical calculus. If one extends the sequent calculus to the capacity of expression of our calculus, then the sequent $a_1,\dotsc,a_m\Imp b_1,\dotsc,b_n$ holds exactly if $a_1\land\dotsb\land a_m\IMP b_1\lor\dotsb\lor b_n$ holds.

(B)~For each formula pair $a,b$ holds: $1\IMP 0\imp a\IMP b$. Thus, if $1\IMP 0$ were to hold, then in the logical interpretation each proposition would imply each other proposition, i.e.\ the calculus would be contradictory. The proof of freedom from contradiction has thus to show that $1\IMP 0$ does not hold.

(C)~The axiom~\pv{(2b)}{(1b)} vs.~\pv{(2c)}{(1c)} corresponds to the classical axiom of comprehension vs.\ extensionality. The axiom~\pv{(2d)}{(1d)} corresponds to the axiom of infinity and postulates the existence of an irreflexive, transitive binary relation in the individual domain whose domain is the whole domain. For our proof of freedom from contradiction, it is of no concern whether this or another equivalent form of the axiom of infinity is postulated.

(D)~The concept of formula is defined constructively. The set of formulae is the smallest set of signs that satisfies the conditions~(F1)--(F3).

Therefore the following ``formula induction'' holds: if a claim holds
\begin{compactenum}[1.]
\item for each prime formula,
\item for $a\land b$, $a\lor b$, $\neg a$, $\forall_{x^\tau}c(x^\tau)$, $\exists_{x^\tau}c(x^\tau)$ if for $a$,~$b$, and~$c(p^\tau)$,
\end{compactenum}
\noindent then it holds for each formula.\bruch

The concept of theorem is defined constructively as well\deutsch{ebenso}. The set of theorems is the smallest set of formula pairs that satisfies conditions~(1) and~(2).

The following theorem induction holds: if a claim holds
\begin{compactenum}[1.]
\item for each axiom,
\item for each conclusion of a rule whose premisses are theorems if for the premisses,
\end{compactenum}
\noindent then it holds for each theorem.

(E)~The following ``duality principle'' holds for the rules of the calculus: if one swaps in each formula pair the left formula with the right one and, in doing so\deutsch{dabei}, simultaneously~$\land$ with~$\lor$ and~$\forall$ with~$\exists$, then each rule transforms again into a rule.

At this\deutsch{hierbei}, negations may remain unchanged. But if one wants to extend the duality also to the theorems, then in addition\deutsch{noch} $a$~must always be swapped with~$\neg a$ and 0~with~1.

(F)~It suffices to restrict oneself to ``proper'' formulae, i.e.\ to formulae that contain neither~0 nor~1. If $c$~is proper, then we write~$\IMP c$ instead of~$1\IMP c$ and $c\IMP$ instead of~$c\IMP0$. We write~$\IMP$ instead of~$1\IMP 0$. If in addition\deutsch{noch} we leave aside the theorems~$0\IMP c$ and~$c\IMP1$, then we obtain a calculus in which only proper formulae occur\deutsch{vorkommen}. We call this calculus the ``proper deductive calculus''.

\section{An inductive calculus.}
We use the same ``signs'' as for the deductive calculus and add~$+$ and~$>$.

We take over the definition of ``types'' and ``free'' vs.\ ``bound variables'' from the deductive calculus.

\begin{sloppypar}
For the ``inductive calculus'' to be constructed we define the concept of constants simultaneously with the concept of formula. As ``constants of type~0'' we take:
\begin{compactenum}[{[}C1{]}]
\item 1.
\item With~$a$ also~$a+1$.
\end{compactenum}
(For communicating the constants of type~0 we use $a^0$, $b^0$, $\dotsc$.) As ``constant of type~$\tau$'' we take for~$\tau=\nu\,(\tau_1\dotsm\tau_n)$ and each formula of $\nu$th order $c(p_1^{\tau_1},\dotsc,p_n^{\tau_n})$ in which $x_1^{\tau_1},\dotsc,x_n^{\tau_n}$ do not occur\deutsch{vorkommen} and which is without variables and constants of types of order~$\nu$: $(x_1^{\tau_1}\dotsm x_n^{\tau_{\pv{1}{n}}})^\nu\,c(x_1^{\tau_1},\dotsc,x_n^{\tau_n})$. (For communicating these constants we use $a^\tau$, $b^\tau$, $\dotsc$.)
\end{sloppypar}

As ``formulae of $\nu$th order'' ($\nu>0$) we take:
\begin{linedisplay}
  [F1]&$(p^0>q^0)$.&\cr
  &$(p^\tau p_1^{\tau_1}\dotsm p_n^{\tau_n})$ for $\tau=\mu(\tau_1\dotsm\tau_n)$ and $\mu\leqq\nu$.&
\end{linedisplay}
\noindent With~$c(p^\tau)$ also $c(a^\tau)$~for each constant~$a^\tau$. These formulae are called the ``prime formulae.'' The formulae~$\pv{p}{a}^0>\pv{q}{b}^0$ are called ``numerical formulae.''
\Item{[F2]}With $a$, $b$, $c(p^\tau)$ also $a\land b$, $a\lor b$, $\neg a$, $\forall_{x^\tau}c(x^\tau)$, $\exists_{x^\tau}c(x^\tau)$ if $x^\tau$~does not occur\deutsch{vorkommen} in~$c(p^\tau)$ and if for the order~$\mu$ of~$\tau$ holds $\mu<\nu$.

The logical interpretation of the formulae is to be carried out as in the deductive calculus. The constants of 0th~type are to be interpreted as the natural numbers. $p^0>q^0$ is to be read as ``$p^0$ greater than $q^0$.'' The constants $a^\tau=(x_1^{\tau_1}\dotsm x_n^{\tau_n})^\nu\,a(x_1^{\tau_1},\dotsc,x_n^{\tau_n})$ are to be interpreted as ``the relation of $\nu$th~order between $x_1^{\tau_1},\dotsc,x_n^{\tau_n}$ defined by $a(x_1^{\tau_1},\dotsc,x_n^{\tau_n})$.''\bruch

We define the concept of theorem as in the deductive calculus for proper formulae:
\Item{[1]}For numerical formulae~$c=(a^0>b^0)$ that are correct vs.\ false on the basis\deutsch{aufgrund} of the interpretation in terms of content of~$a^0$ and~$b^0$ as natural numbers: $\IMP c$ vs.\ $c\IMP$.

\noindent These theorems are called the ``numerical theorems.''
\Item{[2]}\label{ind[2]}Structure rule. If $a_1\land\dotsb\land a_m\IMP b_1\lor\dotsb\lor b_n$ holds, then by the following changes to the left or right formula arises again a theorem: association, i.e.\ the grouping of terms by brackets may be changed; transposition, i.e.\ the sequential arrangement of the terms may be changed; contraction, i.e.\ of several equal terms one may be omitted.
\begin{linedisplay}
  [3a]&$a\IMP b\imp a\IMP b\lor c$.&\cr
  &$a\IMP b\imp a\land c\IMP b$.&\cr
  [3b]&$a\IMP b_1\lor c\comma a\IMP b_2\lor c\imp a\IMP(b_1\land b_2)\lor c$.&\cr
  &$a_1\land c\IMP b\comma a_2\land c\IMP b\imp(a_1\lor a_2)\land c\IMP b$.&\cr
  [3c]&$a\land c\IMP b\imp a\IMP b\lor\neg c$.&\cr
  &$a\IMP b\lor c\imp a\land\neg c\IMP b$.&
\end{linedisplay}
\Item{[3d]}\label{3d}For formulae~$c(p^\tau)$ in which $x^\tau$ does not occur\deutsch{vorkommen}, if $p^\tau$ does not occur\deutsch{vorkommen} in~$a$ nor~$b$:
\[\begin{gathered}
  a\IMP b\lor c(p^\tau)\imp a\IMP b\lor\forall_{x^\tau}c(x^\tau)\text,\\
  a\land c(p^\tau)\IMP b\imp a\land\exists_{x^\tau}c(x^\tau)\IMP b\text.
\end{gathered}\]
For constants~$a^\tau$:
\[\begin{gathered}
   a\land c(a^\tau)\IMP b\imp a\land\forall_{x^\tau}c(x^\tau)\IMP b\text,\\
  a\IMP b\lor c(a^\tau)\imp a\IMP b\lor\exists_{x^\tau}c(x^\tau)\text. 
\end{gathered}\]
In these rules~[3], $a$~and~$b$ may also be omitted if $a\land c$ and $b\lor c$ are replaced by~$c$.
\begin{compactenum}[{[}1{]}]\setcounter{enumi}{3}
\item \label{inductionrule}Induction rule\label{ali:induction}: a formula pair~$C(p^\tau)$ is a theorem
  if for each constant~$a^\tau$ of type~$\pv{c}{\tau}$ the formula
  pair~$C(a^\tau)$ is a theorem.
\item Rule of constants. In order to formulate this rule we need
  the concept of ``elimination of constants.'' Let $a^\tau$ be a constant
  of $\nu$th order
  \[a^\tau=(x_1^{\tau_1}\dotsm x_n^{\tau_n})^\nu\,a(x_1^{\tau_1},\dotsc,x_n^{\tau_n})\text.\]
  By ``elimination'' of~$a^\tau$, we understand the mapping of the set
  of formulae into itself that will be defined in the following
  way. Let $c/a^\tau$ be the image of~$c$.\\
  (i) For each prime formula beginning with~$a^\tau$:
  \[a^\tau c_1\dotsm c_n/a^\tau=a(c_1,\dotsc,c_n)\text.\]
  For each other prime formula~$c$: $c/a^\tau=c$.
  \begin{gather*}
    (c_1\land c_2)/a^\tau=c_1/a^\tau\land c_2/a^\tau\text.\tag{ii}\\
    (c_1\lor c_2)/a^\tau=c_1/a^\tau\lor c_2/a^\tau\text.\\
    \neg c/a^\tau=\neg{c/a^\tau}\text{.\bruch}
  \end{gather*}
  For $c(p^\sigma)/a^\tau=c'(p^\sigma)$:
  \[\begin{gathered}
    (\forall_{x^\sigma}c(x^\sigma))/a^\tau=\forall_{x^\sigma}c'(x^\sigma)\text,\\
    (\exists_{x^\sigma}c(x^\sigma))/a^\tau=\exists_{x^\sigma}c'(x^\sigma)\text.
  \end{gathered}\]
  The image~$c/a^\tau$ is obviously defined by this\deutsch{hierdurch} for each formula~$c$. If $C$~designates the formula pair~$c_1,c_2$, then let $C/a^\tau$ be the formula pair~$c_1/a^\tau,c_2/a^\tau$. Now we can formulate the rule of constants: ``if $C/a^\tau$ is a theorem, then so is~$C$.''
\end{compactenum}
We call [2]--[5] the ``rules'' of the calculus.

The calculus is described by the complete definition of the concept of formula and theorem. We call this calculus shortly the ``inductive calculus.''

Comments to the inductive calculus. [A]~In the inductive calculus holds according to~[3a]
\[{\IMP}\imp a\IMP b\]
for each formula pair $a,b$. The inductive calculus is obviously free from contradiction in the sense that $\IMP$ is not a theorem. In fact, there is no rule that could have $\IMP$ as conclusion. For each rule---except the structure rules---the conclusion contains at least one proper formula. The structure rules trivially cannot have~$\IMP$ as conclusion, as long as the premiss is different from~$\IMP$.

[B]~Instead of the axioms of the deductive calculus appear alone the numerical theorems of the inductive calculus. At\deutsch{bei} defining these theorems, use is made of the interpretation in terms of content.

[C]~The induction rule yields a conclusion out of infinitely many premisses. But the infinite set of premisses is constructively defined, as the set of the constants~$a^\tau$ is defined constructively. Therefore the induction rule is constructively admissible.

[D]~The concept of formula and theorem is again defined constructively as in the deductive calculus.

Therefore ``formula induction of 1st~kind'' also holds: if a claim holds
\begin{compactenum}[1.]
\item for each prime formula,
\item for $a\land b$, $a\lor b$, $\neg a$, $\pv{\exists}{\forall}_{x^\tau}c(x^\tau)$, $\exists_{x^\tau}c(x^\tau)$ if for $a$, $b$, $c(a^\tau)$,
\end{compactenum}
then it holds for each formula.

``Theorem induction'': if a claim holds
\begin{compactenum}[1.]
\item for each numerical theorem,
\item for the conclusion of a rule whose premisses are theorems if for these premisses,
\end{compactenum}
then it holds for each theorem.

It is essential for the inductive calculus that also the following ``formula induction of 2nd~kind'' is valid: if a claim holds
\begin{compactenum}[1.\hfill][10]
\item for numerical formulae,
\item[2.1]for $a\land b$, $a\lor b$, $\neg a$ if for $a$, $b$,
\item[2.2]for~$c$ if for~$c/a^\tau$,
\item[2.3]for $c(p^\tau)$, $\forall_{x^\tau}c(x^\tau)$, $\exists_{x^\tau}c(x^\tau)$ if for each $c(a^\tau)$,
\end{compactenum}
then it holds for each formula.

Proof. A claim that fulfils~1.,~2.,\ holds at first for each formula~($p^0>q^0$) according to~1.\ and~2.3, thus\deutsch{also} for each formula of 1st~order without variables and\bruch{} constants of 1st~order. If the claim holds for each formula of $\nu$th~order without variables and constants of $\nu$th~order, then it holds according to~2.2 for each formula $a^\tau p_1^{\tau_1}\dotsm p_n^{\tau_n}$ in which $a^\tau$~is a constant of order~$\nu$, thus\deutsch{also} according to~2.3 for each prime formula of $\nu$th order, i.e.\ for the prime formulae of $\nu+1$st~order without variables and constants of $\nu+1$st~order. According to~2.1 and~2.3, the validity for each formula of $\nu+1$st~order without variables and constants of $\nu+1$st~order follows from this. Thereby the formula induction of 2nd~kind is proved.

[E]~The same duality principle holds for the rules of the inductive calculus as in the deductive calculus.

\section{The freedom from contradiction of the deductive calculus.}\label{sec:consistency}
We prove that the proper deductive calculus is a part of the inductive calculus. First each proper formula of the deductive calculus is obviously also a formula of the inductive calculus. We extend therefore the proper deductive calculus if we replace its definition of formulae by the definition of formulae of the inductive calculus, but keep its definition of the concept of theorem.

We have then to prove in addition\deutsch{noch} that each theorem of the proper deductive calculus is also a theorem of the inductive calculus. This claim on all deductive theorems is to be proved by a theorem induction. Thus\deutsch{also} we have to prove that
\begin{compactenum}[(I)\hfill][2]
\item the axioms of the proper deductive calculus are inductive theorems,
\item the conclusion of a deductive rule is an inductive theorem if the premisses are inductive theorems.
\end{compactenum}
We shall prove the claims~(I) and~(II) by the formula and theorem inductions valid for the inductive calculus.

(1a) Axiom~$c\IMP c$. We prove by formula induction of 2nd~kind that for each formula~$c$ holds~$c\IMP c$.

1. $c\IMP c$ holds for numerical formulae~$c$, for $\IMP c$ or $c\IMP$ holds, from which in each case arises $c\IMP c$ according to~[3a].

2.1 Let $c_1\IMP c_1$ and $c_2\IMP c_2$ hold. Then follows, because of $c_1\IMP c_1\imp c_1\land c_2\IMP c_1$, \ $c_2\IMP c_2\imp c_1\land c_2\IMP c_2$, and $c_1\land c_2\IMP c_1\comma c_1\land c_2\IMP c_2\imp c_1\land c_2\IMP c_1\land c_2$, also $c_1\land c_2\IMP c_1\land c_2$. As well\deutsch{ebenso} follows $c_1\lor c_2\IMP c_1\lor c_2$. Because of $c\IMP c\imp c\land\neg c\IMP$ and ${c\land\neg c\IMP}\imp\neg c\IMP\neg c$ follows also $\neg c\IMP\neg c$.

2.2 Let $c/a^\tau\IMP c/a^\tau$ hold. Then $c\IMP c$ holds also according to the rule of constants.

2.3 Let $c(a^\tau)\IMP c(a^\tau)$ hold for each~$a^\tau$. Then follows first  $c(p^\tau)\IMP c(p^\tau)$ according to the induction rule. Further follows, because of $c(a^\tau)\IMP c(a^\tau)\imp\forall_{x^\tau}c(x^\tau)\IMP c(a^\tau)$, $\forall_{x^\tau}c(x^\tau)\IMP c(a^\tau)$ for each~$a^\tau$; thus\deutsch{also} $\forall_{x^\tau}c(x^\tau)\IMP c(p^\tau)$ according to the induction rule, and from this $\forall_{x^\tau}c(x^\tau)\IMP\forall_{x^\tau}c(x^\tau)$. As well\deutsch{ebenso} follows $\exists_{x^\tau}c(x^\tau)\IMP\exists_{x^\tau}c(x^\tau)$.

\begin{sloppypar}
(1b) Axiom of comprehension: for $\tau=\nu\,(\tau_1\dotsm\tau_n)$ and formulae of $\nu$th order $c=c(p_1^{\tau_1},\allowbreak\dotsc,p_n^{\tau_n})$,
\[\IMP\exists_{x^\tau}\forall_{x_1^{\tau_1}}\dotsm\forall_{x_n^{\tau_n}}(x^\tau x_1^{\tau_1}\dotsm x_n^{\tau_n}\Equ c(x_1^{\tau_1},\dotsc,x_n^{\tau_n}))\text.\]
According to observation~(1a) holds $c(p_1^{\tau_1},\dotsc,p_n^{\tau_n})\IMP c(p_1^{\tau_1},\dotsc,p_n^{\tau_n})$. From this follows $\IMP c(p_1^{\tau_1},\dotsc,p_n^{\tau_n})\Imp c(p_1^{\tau_1},\dotsc,p_n^{\tau_n})$, $\IMP c(p_1^{\tau_1},\dotsc,p_n^{\tau_n})\Equ c(p_1^{\tau_1},\dotsc,p_n^{\tau_n})$, thus\deutsch{also} $\IMP{\forall_{x_1^{\tau_1}}\dotsm\forall_{x_n^{\tau_n}}({c(x_1^{\tau_1},\dotsc,x_n^{\tau_n})\Equ\allowbreak c(x_1^{\tau_1},\dotsc,x_n^{\tau_n})})}$. According to the rule of constants\bruch{}  follows for $c^\tau=(y_1^{\tau_1}\dotsm y_n^{\tau_n})^\nu\,c(y_1^{\tau_1},\dotsc,y_n^{\tau_n})$, if $c(p_1^{\tau_1},\dotsc,p_n^{\tau_n})$ is without variables and constants of $\nu$th order:
\[\IMP\forall_{x_1^{\tau_1}}\dotsm\forall_{x_n^{\tau_n}}(c^\tau x_1^{\tau_1}\dotsm x_n^{\tau_n}\Equ c(x_1^{\tau_1},\dotsc,x_n^{\tau_n}))\text.\]
From this
\[\IMP\exists_{x^\tau}\forall_{x_1^{\tau_1}}\dotsm\forall_{x_n^{\tau_n}}(x^\tau x_1^{\tau_1}\dotsm x_n^{\tau_n}\Equ c(x_1^{\tau_1},\dotsc,x_n^{\tau_n}))\text.\]
If $c$~contains a constant~$a^\sigma$ of order~$\nu$ and the axiom of comprehension holds for $c/a^\sigma$, then according to~[5] also for~$c$. If $c=c(p^\sigma)$ contains a free variable~$p^\sigma$ of order~$\nu$ and the axiom of comprehension holds for all $c(a^\sigma)$, then according to~[4] also for~$c$. From this follows the axiom of comprehension for each formula of $\nu$th order.
\end{sloppypar}

(1c) Axiom of extensionality:
\[\forall_{x_1^{\tau_1}}\dotsm\forall_{x_n^{\tau_n}}(p^\tau x_1^{\tau_1}\dotsm x_n^{\tau_n}\Equ q^\tau x_1^{\tau_1}\dotsm x_n^{\tau_n})\IMP c(p^\tau)\Imp c(q^\tau)\text.\]
We write for the left formula for abbreviating $p^\tau\equiv q^\tau$, and shall for each~$c(p^\tau)$ prove $p^\tau\equiv q^\tau\IMP c(p^\tau)\Equ c(q^\tau)$. We use for this the following induction, that results immediately from the formula induction of 2nd~kind: if a claim holds\\
1. for each prime formula $p^\tau a_1^{\tau_1}\dotsm a_n^{\tau_n}$ and for each formula that does not contain~$p^\tau$,\\
2.1 for $a\land b$, $a\lor b$, $\neg a$ if for $a$, $b$,\\
2.2 for~$c$ if for~$c/b^\sigma$,\\
2.3 for $c(q^\sigma)$, $\forall_{x^\sigma}c(x^\sigma)$, $\exists_{x^\sigma}c(x^\sigma)$ if for each $c(b^\sigma)$ and $b^\sigma\ne p^\tau$,\\
then it holds for each formula~$c$.

1. According to~(1a) holds
\[p^\tau a_1^{\tau_1}\dotsm a_n^{\tau_n}\Equ q^\tau a_1^{\tau_1}\dotsm a_n^{\tau_n}\IMP
  p^\tau a_1^{\tau_1}\dotsm a_n^{\tau_n}\Equ q^\tau a_1^{\tau_1}\dotsm a_n^{\tau_n}\text,\]
from which follows $p^\tau\equiv q^\tau\IMP p^\tau a_1^{\tau_1}\dotsm a_n^{\tau_n}\Equ q^\tau a_1^{\tau_1}\dotsm a_n^{\tau_n}$.

For~2.1--2.3 it suffices to prove the following:
\[\begin{gathered}
  a\IMP c_1\Equ c_2\comma a\IMP d_1\Equ d_2\imp a\IMP (c_1\land \pv{c_2}{d_1})\Equ(\pv{d_1}{c_2}\land d_2)\text.\tag{a}\\
  a\IMP c_1\Equ c_2\comma a\IMP d_1\Equ d_2\imp a\IMP (c_1\lor \pv{c_2}{d_1})\Equ(\pv{d_1}{c_2}\lor d_2)\text.\\
  a\IMP c_1\Equ c_2\imp a\IMP \neg c_1\Equ \neg c_2\text.
\end{gathered}\]
(b) $a\IMP \pv{a_1}{c_1}/b^\sigma\Equ c_2/b^\sigma\imp a\IMP c_1\Equ c_2$ if $b^\sigma$~does not occur\deutsch{vorkommen} in~$a$.

\noindent(c) If $a\IMP c_1(b^\sigma)\Equ c_2(b^\sigma)$ for each $b^\sigma$  and $q^\sigma$ does not occur\deutsch{vorkommen} in~$a$, then
\[\begin{gathered}
  a\IMP c_1(q^\sigma)\Equ c_2(q^\sigma)\text,\\
  a\IMP\forall_{x^\sigma}c_1(x^\sigma)\Equ\forall_{x^\sigma}c_2(x^\sigma)\text,\qquad
  a\IMP\exists_{x^\sigma}c_1(x^\sigma)\Equ\exists_{x^\sigma}c_2(x^\sigma)\text.
\end{gathered}\]

We need for these two auxiliary rules.

For $a\IMP c\land d\imp a\IMP c$, we prove the stronger \textsc{auxiliary rule~1}:
\[a\IMP b\lor(c\land d)\lor\dotsb\lor(c\land d)\imp a\IMP b\lor c\text.\]

This is a claim on all inductive theorems. It claims that for each theorem holds: if a theorem~$C$ has the form $a\IMP b\lor(c\land d)\lor\dotsb\lor(c\land d)$, then holds $a\IMP b\lor c$. For the proof we may therefore apply a theorem induction.\bruch

[1] For numerical theorems nothing is to be proved. We consider the inductive rules and assume the validity of the claim for the premisses.

[2] Structure rules. Let $C$~be the conclusion of a structure rule. If one applies the induction hypothesis to the premisses, then one obtains a theorem from which may at once be inferred $a\IMP b\lor c$ by a structure rule.

[3a] If $C$~is conclusion of an inductive rule~[3a], then the premiss has the form $a_1\IMP b_1\lor(c\land d)\lor\dotsb\lor(c\land d)$, and from $a_1\IMP b_1\lor c$ (possibly $a_1\IMP b_1$) may be inferred $a\IMP b\lor c$.

[3b] If $C$~is conclusion of a rule~[3b], then we have as premisses $a\IMP b\lor c\lor(c\pv{\lor}{\land} d)\lor\dotsb$, $a\IMP b\lor d\lor(c\land d)\lor\dotsb$, or yet premisses of the form $a_1\IMP b_1\lor(c\land d)\lor\dotsb$, so that from $a_1\IMP b_1\lor c$ may at once again also $a\IMP b\lor c$ be inferred by a rule~[3b]. In the first case follows $a\IMP b\lor c\lor c$, thus\deutsch{also} $a\IMP b\lor c$.

[3c] If $C$~is conclusion of a rule~[3c], then the premiss has the form $a_1\IMP b_1\lor(c\land d)\lor\dotsb\lor(c\land d)$, and from $a_1\IMP b_1\lor c$ we may infer $a\IMP b\lor c$.

[3d] The same as for~[3c] holds verbatim.

[4] If $C=C(p^\tau)$ and $C$~is conclusion out of the premisses~$C(a^\tau)$ for each~$a^\tau$, then the induction hypothesis yields---if we designate the pair~$a,b\lor c$ by~$D=D(p^\tau)$---at once that $D(a^\tau)$ is a theorem for each~$a^\tau$. Thus\deutsch{also} $D$~is also a theorem.

[5] If $C$~is conclusion of a rule of constants, then we have the premiss~$C/a^\tau$, from which follows according to the induction hypothesis that $D/a^\tau$ is also a theorem. Thus\deutsch{also} $D$~is also a theorem.

Thereby auxiliary rule~1 is proved. As one sees, all of the single steps are trivial. We shall therefore not treat them anymore in the following similar proofs.

For $a\IMP b\lor\neg c\imp a\land c\IMP b$, we prove by theorem induction the stronger \textsc{auxiliary rule~2}:
\[a\IMP b\lor\neg c\lor\dotsb\lor\neg c\imp a\land c\IMP b\text.\]

For numerical theorems nothing is to be proved. The treatment of rules [2]--[5] is trivial in all cases. Let only [3c]~be emphasised:
\[a\land c\IMP b\lor\neg c\lor\dotsb\lor\neg c\imp a\IMP b\lor\neg c\lor\neg c\lor\dotsb\lor\neg c\text.\]
The induction hypothesis yields $a\land c\land c\IMP b$, thus\deutsch{also} $a\land c\IMP b$ also follows. Thereby auxiliary rule~2 is proved. 

Under addition of these auxiliary rules to the rules of the inductive calculus, we may now undertake the formula induction for the axiom of extensionality, and indeed with exactly the same inferences as for the axiom~$c\IMP c$. From $a\IMP c_1\Equ c_2$ we may now in fact first infer $a\IMP c_1\Imp c_2$ and then $a\land c_1\IMP c_2$. Everything remaining is then to be concluded as under~(1a).

(1d) Axiom of infinity: for $\tau=1(00)$,
\[\IMP\exists_{\!\!x^\tau}(\forall_{x^0}\neg{x^\tau x^0x^0}\land\forall_{x^0}\!\!\exists_{\!\!y^0}x^\tau x^0y^0\land\forall_{x^0}\!\forall_{y^0}\!\forall_{z^0}((x^\tau x^0y^0\land x^\tau y^0z^0)\Imp x^\tau x^0z^0))\text{.\bruch}\]
We prove for $a^\tau=(u^0v^0)^1\,v^0>u^0$
\[\begin{gathered}
  \IMP\forall_{x^0}\neg{a^\tau x^0x^0},\qquad\IMP\forall_{x^0}\pv{\forall}{\exists}_{y^0}a^\tau x^0y^0\text,\\
  \IMP\forall_{x^0}\forall_{y^0}\forall_{z^0}((a^\tau x^0y^0\land a^\tau y^0z^0)\Imp a^\tau x^0z^0)\text;
\end{gathered}\]
from which the axiom follows at once. Because of the rule of constants it suffices to prove
\[\begin{gathered}
  \IMP\forall_{x^0}\neg{x^0>x^0},\qquad\IMP\forall_{x^0}\exists_{y^0}y^0>x^0\text,\\
  \IMP\forall_{x^0}\forall_{y^0}\forall_{z^0}((y^0>x^0\land z^0>y^0)\Imp z^0>x^0)\text.
\end{gathered}\]
For each~$a^0$ holds $a^0>a^0\IMP$; from which $\IMP\neg{a^0>a^0}$, and according to the induction rule $\IMP\neg{p^0>p^0}$, thus\deutsch{also} $\IMP\forall_{x^0}\neg{x^0>x^0}$ follows. Moreover\deutsch{ferner} $\IMP a^0+1>a^0$ holds for each~$a^0$; from which $\IMP\exists_{y^0}y^0>a^0$ and according to the induction rule $\IMP\exists_{y^0}y^0>p^0$. Thus\deutsch{also} $\IMP\forall_{x^0}\exists_{y^0}y^0>x^0$ follows. For each $a^0$, $b^0$, $c^0$ holds finally $\IMP c^0>a^0$ or $c^0>b^0\IMP$ or $b^0>a^0\IMP$; from which in each case follows $b^0>a^0\land c^0>b^0\IMP c^0>a^0$ according to~[3a]. Thus\deutsch{also} holds also $q^0>p^0\land r^0>q^0\IMP r^0>p^0$ according to the induction rule, from which \pv{$\IMP(q^0>p^0\land r^0>q^0\Imp r^0>p^0)$}{$\IMP(q^0>p^0\land r^0>q^0)\Imp r^0>p^0$} and $\IMP\forall_{x^0}\forall_{y^0}\forall_{z^0}((y^0>x^0\land z^0>y^0)\Imp z^0>x^0)$ follows.

(2) Now it still remains\deutsch{bleiben} to prove the deductive rules for the inductive calculus.

(2a) $a\IMP c\comma c\IMP b\imp a\IMP b$.

We prove the stronger \textsc{auxiliary rule~3}:
\[a_1\IMP b_1\lor c\lor\dotsb\lor c\comma a_2\land c\land\dotsb\land c\IMP b_2\imp a_1\land a_2\IMP b_1\lor b_2\text.\]

For this we use the formula induction of 2nd~kind for~$c$.

1. Let $c$~be a numerical formula. We prove auxiliary rule~3 by theorem induction for each theorem $a_1\IMP b_1\lor c\lor\dotsb\lor c$.

1.1 Let $a_1\IMP b_1\lor c\lor\dotsb\lor c$ be a numerical theorem. We prove auxiliary rule~3 by theorem induction for each theorem $a_2\land c\land\dotsb\land c\IMP b_2$.

1.1.1 Let $a_2\land c\land\dotsb\land c\IMP b_2$ be a numerical theorem. Auxiliary rule~3 is valid, because $\IMP c$ and $c\IMP$ do not hold simultaneously.

1.1.2 Let $a_2\land c\land\dotsb\land c\IMP b_2$ be a conclusion of an inductive rule, and let auxiliary rule~3 be valid for the premisses. The treatment of rules [2]--[5] is trivial in each case because $c$~is a numerical formula.

1.2 Let $a_1\IMP b_1\lor c\lor\dotsb\lor c$ be conclusion of an inductive rule, and let auxiliary rule~3 be valid for the premisses. The treatment of rules [2]--[5] is trivial as under~1.1.2.

2.1 Let auxiliary rule~3 be valid for $c=c_1$ and $c=c_2$. We prove it from this for $c=c_1\land c_2$, $c=c_1\lor c_2$, and $c=\neg c_1$. According to auxiliary rule~1 holds
\[\begin{gathered}
  a_1\IMP b_1\lor(c_1\land c_2)\lor\dotsb\lor(c_1\land c_2)\imp a_1\IMP b_1\lor c_1\text,\\
  a_1\IMP b_1\lor(c_1\land c_2)\lor\dotsb\lor(c_1\land c_2)\imp a_1\IMP b_1\lor c_2\text.
\end{gathered}\]
Moreover\deutsch{ferner} holds $a_2\land(c_1\land c_2)\land\dotsb\land(c_1\land c_2)\IMP b_2\imp a_2\land c_1\land c_2\IMP b_2$. But according to hypothesis holds $a_1\IMP b_1\lor c_1\comma a_2\land c_1\land c_2\IMP b_2\imp a_1\land a_2\land c_2\IMP b_1\lor b_2$ and $a_1\IMP b_1\lor c_2\mathbin{\text{\sffamily\upshape\bfseries,\bruch}}a_1\land a_2\land c_2\IMP b_1\lor b_2\imp a_1\land a_1\land a_2\IMP b_1\lor b_1\lor b_2$. By contraction arises $a_1\land a_2\IMP b_1\lor b_2$. Dually to auxiliary rule~1 holds $a\land(c_1\lor c_2)\land\dotsb\land(c_1\lor c_2)\IMP b\imp a\land c_1\IMP b$, and by its aid follows as just also the validity of auxiliary rule~3 for $c=c_1\lor c_2$. Finally holds according to auxiliary rule~2, $a_1\IMP b_1\lor\neg c_1\lor\dotsb\lor\neg c_1\imp a_1\land c_1\IMP b_1$; and dually to that holds $a_2\land\neg c_1\land\dotsb\land\neg c_1\IMP b_2\imp a_2\IMP b_2\lor c_1$. But according to hypothesis holds $a_2\IMP b_2\lor c_1\comma a_1\land c_1\IMP b_1\imp a_1\land a_2\IMP b_1\lor b_2$.

2.2 Let auxiliary rule~3 be valid for $c=d/a^\tau$. We prove it from this for $c=d$.

For this we use \textsc{auxiliary rule~4}: if $C$~is a theorem, then also $C/a^\tau$.

The proof by theorem induction is trivial for each step.

We designate the theorems appearing\deutsch{auftretend} in auxiliary rule~3 by $C_1$, $C_2$, and~$C_3$, so that it states: if $C_1$~and~$C_2$ are theorems, then also~$C_3$. According to auxiliary rule~4 follows first that $C_1/a^\tau$~and~$C_2/a^\tau$ are theorems. According to hypothesis follows from that that $C_3/a^\tau$~is a theorem. According to the rule of constants, then also $C_3$ is a theorem.

2.3 Let auxiliary rule~3 be valid for each~$c=d(a^\tau)$, where $a^\tau$ runs through\deutsch{durchläuft} all constants of type~$\tau$.

We prove it from this for~$c=d(p^\tau)$ exactly correspondingly to~2.2 by using \textsc{auxiliary rule~5}: if $C(p^\tau)$~is a theorem, then also $C(a^\tau)$ for each~$a^\tau$.

The proof is again trivial for each step.

It remains\deutsch{bleibt} yet to show the validity of auxiliary rule~3 for $c=\forall_{x^\tau}d(x^\tau)$ and $c=\exists_{x^\tau}d(x^\tau)$.

First, we prove again by a trivial theorem induction \textsc{auxiliary rule~6}:
\[a_1\IMP b_1\lor\forall_{x^\tau}d(x^\tau)\lor\dotsb\lor\forall_{x^\tau}d(x^\tau)\imp a_1\IMP b_1\lor d(p^\tau)\text.\]

If we choose for $p^\tau$ a variable that does not occur\deutsch{vorkommen} in $a_1$, $b_1\lor\forall_{x^\tau}d(x^\tau)$, then auxiliary rule~5 yields $a_1\IMP b_1\lor d(a^\tau)$ for each~$a^\tau$. We have to show now---under hypothesis of auxiliary rule~3 for $c=d(a^\tau)$ and the validity of $a_1\IMP b_1\lor d(a^\tau)$ for each $a^\tau$---:
\[a_2\land\forall_{x^\tau}d(x^\tau)\land\dotsb\land\forall_{x^\tau}d(x^\tau)\IMP b_2\imp a_1\land a_2\IMP b_1\lor b_2\text.\]
We prove this claim by theorem induction.

2.3.1 Let $a_2\land\forall_{x^\tau}d(x^\tau)\land\dotsb\IMP b_2$ be a numerical theorem. The claim is then trivially valid.

2.3.2 Let $a_2\land\forall_{x^\tau}d(x^\tau)\land\dotsb\IMP b_2$ be conclusion of an inductive rule, and let the claim be valid for the premisses. The treatment of each rule is trivial, except the one case of rule~[3d]:
\[a_2\land d(a^\tau)\land\forall_{x^\tau}d(x^\tau)\land\dotsb\IMP b_2\imp a_2\land\forall_{x^\tau}d(x^\tau)\land\dotsb\IMP b_2\text.\]
As the claim is valid for the premisses, $a_1\land a_2\land d(a^\tau)\IMP b_1\lor b_2$ follows. Further holds according to hypothesis $a_1\IMP b_1\lor d(a^\tau)$ and
\[a_1\IMP b_1\lor d(a^\tau)\comma a_1\land a_2\land d(a^\tau)\IMP b_1\lor b_2\imp a_1\land a_1\land a_2\IMP b_1\lor b_1\lor b_2\text,\]
from which $a_1\land a_2\IMP b_1\lor b_2$ follows.

Thereby auxiliary rule~3 is proved for $c=\forall_{x^\tau}d(x^\tau)$.

The proof for $c=\exists_{x^\tau}d(x^\tau)$ proceeds dually to this; and auxiliary rule~3 is proved in general.\bruch

The further rules of the deductive calculus now make no difficulties anymore. The deductive rules~(2b) and~(2c) are contained in the inductive rules~[2], [3a]--[3c]. Of the deductive rules~(2d) two are contained in the inductive rules~[3d], the other two in auxiliary rule~6 and the dual auxiliary rule.

Thereby all deductive axioms are recognised as inductive theorems, and all deductive rules as also valid in the inductive calculus.

Thus\deutsch{also} each proper deductive theorem is also an inductive theorem; in particular $\IMP$~is not a theorem in the proper deductive calculus, because $\IMP$ is not a theorem in the inductive calculus. Thereby the freedom from contradiction of the deductive calculus is proved.\phantomsection\label{sec:consistency:end}

\section{The independence of the axiom of reducibility.}
With the method of~\S\pv37 one may also prove the freedom from contradiction of other similar calculuses. The relation of identity may e.g.\ be added to the deductive calculus, and the axiom of extensionality replaced by the axioms of identity:\\
$(c_1)$ $\IMP p^\tau=p^\tau$.\\
$(c_2)$ $p^\tau=q^\tau\IMP c(p^\tau)\Imp c(q^\tau)$.\\
For the proof of freedom from contradiction one has then to modify the inductive calculus in the following way. To the prime formulae one adds~$p^\tau=q^\tau$, to the numerical formulae $a^\tau=b^\tau$. These numerical formulae are interpreted as ``$a^\tau$ and $b^\tau$ are signs of the same form.'' The axiom~$(c_1)$ is then an inductive theorem, because $\IMP a^\tau=a^\tau$ is a numerical theorem for each constant~$a$.

$(c_2)$ follows as well\deutsch{ebenso}, as for two constants~$a^\tau$, $b^\tau$ not of the same the form always holds $a^\tau=b^\tau\IMP$; but for constants~$a^\tau$, $b^\tau$ of the same form $c(a^\tau)\IMP c(b^\tau)$, thus\deutsch{also} $\IMP c(a^\tau)\Imp c(b^\tau)$. According to~[3a] follows in each case $a^\tau=b^\tau\IMP c(a^\tau)\Imp c(b^\tau)$.

The freedom from contradiction of the deductive calculus is even conserved if one adds axioms that postulate in terms of content the equipotence of the set of individuals with the set of relations of type~$\tau$:\\
(e) For $\sigma={\nu+1\,(\tau~0)}$, if $\tau$~of order~$\nu$,
\[\IMP\exists_{x^\sigma}(\forall_{x^\tau}\exists_{x^0}x^\sigma x^\tau x^0\land\forall_{z^0}\forall_{x^\tau}\forall_{y^\tau}((x^\sigma x^\tau z^0\land x^\sigma y^\tau z^0)\Imp x^\tau=y^\tau))\text.\]
For the proof of freedom from contradiction we modify the inductive calculus once more. To the prime formulae $p^\tau\rel p^0$ is being added, to the numerical formulae $a^\tau\rel a^0$. In order to interpret these numerical formulae, we carry out an enumeration of the constants of type~$\tau$. Such an enumeration is possible, as the set of all formulae is countable. We interpret~$a^\tau\rel a^0$ then as ``$a^\tau$ has in the enumeration of the constants of type~$\tau$ the number~$a^0$.'' On the basis\deutsch{aufgrund} of this interpretation there is to each~$a^\tau$ an~$a^0$ with~$a^\tau\rel a^0$. Thus\deutsch{also} for each~$a^\tau$ holds $\IMP\exists_{x^0}a^\tau\rel x^0$; from which follows $\IMP\exists_{x^0}p^\tau\rel x^0$ and $\IMP\forall_{x^\tau}\exists_{x^0}x^\tau\rel x^0$. For arbitrary constants $a^\tau$, $b^\tau$, and $c^0$ holds moreover\deutsch{ferner} $\IMP a^\tau=b^\tau$ or $a^\tau\rel c^0\IMP$ or $b^\tau\rel c^0\IMP$. In each case results according to~[3a] $a^\tau\rel c^0\land b^\tau\rel c^0\IMP a^\tau=b^\tau$, from which follows
\[\begin{gathered}
  p^\tau\rel r^0\land q^\tau\rel r^0\IMP p^\tau=q^\tau\text,\qquad\IMP(p^\tau\rel r^0\land q^\tau\rel r^0)\Imp p^\tau=q^\tau\text,\\
  \IMP\forall_{\pv{x}{z}^0}\forall_{x^\tau}\forall_{y^\tau}((x^\tau\rel z^0\land y^\tau\rel z^0)\Imp x^\tau=y^\tau)\text{.\bruch}
\end{gathered}\]
According to the rule of constants holds therefore for~$a^\sigma=(u^\tau v^0)^{\nu+1}\,u^\tau\rel v^0$,
\[\IMP\forall_{x^\tau}\exists_{x^0}a^\sigma x^\tau x^0\text,\qquad\forall_{z^0}\forall_{x^\tau}\forall_{y^\tau}((a^\sigma x^\tau z^0\land a^\sigma y^\tau z^0)\Imp x^\tau=y^\tau)\text.\]
If we bind these two formulae by~$\land$, then a formula $d(a^\sigma)$ arises; and $\IMP d(a^\sigma)$ holds, from which $\IMP\exists_{x^\sigma}d(x^\sigma)$ follows. This is axiom~(e).

We show in the end that in the modified deductive calculus thereby proven to be free from contradiction, the axiom of reducibility is refutable. A simple case of this axiom is: for~$\rho={\nu+2\,(0)}$, if $\tau=\nu\,(0)$, $\IMP\forall_{x^\rho}\exists_{x^\tau}\pv{\exists}{\forall}_{x^0}x^\rho x^0\Equ x^\tau x^0$. First holds for $\sigma={\nu+1\,(\tau~0)}$, $d(p^\sigma)\IMP(p^\sigma p^\tau r^0\land p^\sigma q^\tau r^0)\Imp p^\tau=q^\tau$ and $p^\tau=q^\tau\IMP\neg{q^\tau r^0}\Imp\neg{p^\tau r^0}$; from which follows $d(p^\sigma)\land p^\sigma p^\tau r^0\land p^\sigma q^\tau r^0\land\neg{q^\tau r^0}\IMP\neg{p^\tau r^0}$ and $d(p^\sigma)\land p^\sigma p^\tau r^0\land\exists_{x^\tau}{(p^\sigma x^\tau r^0\land\neg{x^\tau r^0})}\IMP\neg{p^\tau r^0}$. Because of $p^\sigma p^\tau r^0\land\neg{p^\tau r^0}\IMP\exists_{x^\tau}(p^\sigma x^\tau r^0\land\neg{x^\tau r^0})$ results elementarily
\[\begin{gathered}
  d(p^\sigma)\land p^\sigma p^\tau r^0\land(p^\tau r^0\Equ\exists_{x^\tau}(p^\sigma x^\tau r^0\land\neg{x^\tau r^0}))\IMP\text,\\
  d(p^\sigma)\land p^\sigma p^\tau r^0\land(p^\rho r^0\Equ\exists_{x^\tau}(p^\sigma x^\tau r^0\land\neg{x^\tau r^0}))\land(p^\rho r^0\Equ p^\tau r^0)\IMP\text.
\end{gathered}\]
Further follows now
\[d(p^\sigma)\land\exists_{x^0}p^\sigma p^\tau x^0\land\forall_{x^0}(p^\rho x^0\Equ\exists_{x^\tau}(p^\sigma x^\tau x^0\land\neg{x^\tau x^0}))\land\forall_{x^0}(p^\rho x^0\Equ p^\tau x^0)\IMP\text;\]
thus\deutsch{also} because of~$d(p^\sigma)\IMP\forall_{x^\tau}\exists_{x^0}p^\sigma x^\tau x^0$,
\[d(p^\sigma)\land\forall_{x^0}(p^\rho x^0\Equ\exists_{x^\tau}(p^\sigma x^\tau x^0\land\neg{x^\tau x^0}))\land\exists_{x^\tau}\forall_{x^0}(p^\rho x^0\Equ x^\tau x^0)\IMP\text;\]
because of $\pv{}{\IMP}\exists_{x^\rho}\forall_{x^0}(x^\rho x^0\Equ\exists_{x^\tau}(p^\sigma x^\tau x^0\land\neg{x^\tau x^0}))$,
\[d(p^\sigma)\land\forall_{x^\rho}\exists_{x^\tau}\forall_{x^0}(x^\rho x^0\Equ x^\tau x^0)\pv{}{\IMP}\text;\]
and because of $\IMP\exists_{x^\sigma}d(x^\sigma)$,
\[\forall_{x^\rho}\exists_{x^\tau}\forall_{x^0}(x^\rho x^0\Equ x^\tau x^0)\IMP\text,\]
q.e.d.\medskip

\textsc{university of bonn}
\end{document}